\documentclass[11pt, oneside]{article}
\usepackage{amsfonts}
\usepackage{mathrsfs}
\usepackage{color}
\usepackage[colorlinks]{hyperref}
\usepackage{latexsym}
\usepackage{amssymb}
\usepackage{amsmath}
\usepackage{enumerate}
\usepackage{amsthm}
\usepackage{indentfirst}
\usepackage{mathtools}
\usepackage{tikz}
\usepackage{psfrag}
\usepackage{graphicx}
\usepackage{bm}
\usepackage[rflt]{floatflt}
\usepackage{float}
\usepackage{authblk}
\usepackage[OT2,T1]{fontenc}
\usepackage{ulem}
\usepackage[square, comma, sort&compress, numbers]{natbib}
\usepackage{verbatim}
\newtheorem{problem}{Problem}

\numberwithin{equation}{section}
\allowdisplaybreaks

\newtheorem{theorem}{\color{black}\indent Theorem}[section]
\newtheorem{lemma}{\color{black}\indent Lemma}[section]
\newtheorem{proposition}{\color{black}\indent Proposition}[section]

\newtheorem{remark}{\color{black}\indent Remark}[section]
\newtheorem{corollary}{\color{black}\indent Corollary}[section]

\usepackage{amssymb,amsmath}
\begin{document}
\title{\LARGE\bf Global Existence of Classical Solutions to the
Relativistic Quantum Hydrodynamic System with
Small Initial Data}

\author[1]{Ben Duan}

\author[1,2]{Rongrong Yan\thanks{Corresponding author.
Email: rryan.math@gmail.com}}

\affil[1]{\small School of Mathematics, Jilin University,
Qianjin Road, Changchun 130012, Jilin, People's Republic of China}

\affil[2]{\small Department of Mathematical Sciences,
Seoul National University, Seoul 08826, Republic of Korea}

\date{}

\maketitle

\begin{abstract}
We establish global existence, decay, and scattering for sufficiently
small, smooth, and localized perturbations of a constant non-vacuum
equilibrium of a relativistic quantum hydrodynamic system in three
space dimensions. In logarithmic-amplitude and phase variables,
the equations form a semilinear system of coupled wave equations.
The skew-symmetric coupling between the time derivatives cancels
in the energy identity, yielding a natural derivative energy. The linearized system has two dispersion branches. At low frequencies,
the slow branch exhibits Schr\"odinger-type dispersion, while the fast
branch has a spectral gap; at high frequencies, both branches are
wave-like. The main nonlinear difficulty arises from quadratic
interactions with nontrivial time and space-time resonances.
We show that the symbol of every active quadratic interaction
contains the corresponding interaction phase as an exact factor.
This structural cancellation removes the resonant denominator
and allows us to eliminate the quadratic terms by a nonsingular
normal-form transformation. Combining this transformation with
dispersive and energy estimates, we obtain uniform-in-time bounds
for the derivative energy and $\langle t\rangle^{-3/2}$ decay of
the first derivatives in $L^\infty$. We further prove that the
nonlinear solution scatters to a solution of the linearized
system in the high-order derivative energy norm.

\end{abstract}



\section{Introduction}
Lin and Wu \cite{LinWu2012} introduced this hydrodynamic model for irrotational flows in their study of the hydrodynamic limit of the modulated nonlinear Klein--Gordon equation. In this setting, the velocity field is given by \(u=\nabla S\), and the relativistic quantum hydrodynamic system (RQHD) \cite{DuanYan} can be equivalently expressed in terms of the particle density \(n\) and the phase \(S\) as follows:
\begin{align}\label{J1}
\begin{cases}
\partial_t n+\operatorname{div}(n\nabla S)
=\upsilon^2\partial_t(n\partial_t S),\\[2mm]
\partial_t(n\nabla S)
+\operatorname{div}\!\left(
\dfrac{n\nabla S\otimes n\nabla S}{n}
\right)
-\dfrac{\varepsilon^2}{2}n\nabla\!\left(
\dfrac{\Delta\sqrt n}{\sqrt n}
\right)
\\[2mm]
\qquad\qquad\qquad\qquad\qquad\qquad
=\frac{1}{2}\upsilon^2\left[
\partial_t\bigl(2n\partial_tS\nabla S\bigr)
-\varepsilon^2n\nabla\!\left(
\dfrac{\partial_{tt}\sqrt n}{\sqrt n}
\right)
\right].
\end{cases}
\end{align}
Here $n$ and $S$ denote the particle density and the phase, respectively, while $\varepsilon>0$ and $\upsilon>0$ denote the quantum and relativistic parameters. Although the explicit parameters make the non-relativistic and semiclassical reductions transparent, our focus is the long-time dynamics at fixed parameter values. We therefore introduce the logarithmic density $\psi=\log\sqrt n$, $n=e^{2\psi}$,
and, after a normalization of the fixed parameters, obtain the system \eqref{eq:main-system}. Under this reformulation, the normalized non-vacuum state \(n=1\) corresponds to \(\psi=0\). The central problem is then to determine whether small, smooth, and localized perturbations of \((\psi,S)=(0,0)\) can be controlled globally by the dispersive effects of the pressureless relativistic quantum system. Such global control keeps the density uniformly away from vacuum and leads to decay and scattering of the physical variables.

Earlier studies addressed the derivation of the hydrodynamic equations and the analysis of the associated parameter-dependent singular limits \cite{LinWu2012}, while the local existence and uniqueness of non-vacuum classical solutions have only recently been established for several related models, including Poisson-coupled systems \cite{DuanYan}. These results, however, do not address the global-in-time dynamics of the pressureless Euclidean system considered in this paper.

The mathematical analysis of quantum fluid models has attracted considerable attention over the past decades. For the inviscid Euler--Korteweg system, the local well-posedness of strong solutions was established by Benzoni-Gavage, Danchin and Descombes \cite{BenzoniGavageDanchinDescombes2007}. At the finite-energy level, Antonelli and Marcati constructed global weak solutions for quantum hydrodynamic systems \cite{Antonelli2009} and subsequently developed a two-dimensional theory that remains meaningful in the presence of vacuum regions \cite{Antonelli2012}. Related global weak-solution theories have also been established for dissipative quantum Navier--Stokes equations \cite{Antonelli2017} and for quantum fluids coupled to electromagnetic fields \cite{Antonelli2024}. These works provide a broad existence theory for quantum fluid models, including large-data and vacuum regimes, but they do not address the global smooth dynamics, decay, or scattering of small localized perturbations.

The global behavior of small solutions to nonlinear hyperbolic and dispersive equations has been extensively studied. For nonlinear wave equations in three space dimensions, the vector-field method and the null condition developed in the classical works of Klainerman and Christodoulou provide fundamental mechanisms for the small-data global theory; see, for instance, \cite{Klainerman1985,Christodoulou1986}. For quadratic Klein--Gordon equations, Shatah introduced the normal-form method \cite{Shatah1985}, which exploits temporal oscillations to remove suitable quadratic interactions. The role of resonances in nonlinear dispersive equations was later systematically organized through the space--time resonance method; see, for example, \cite{GermainMasmoudiShatah2012,GermainMasmoudiShatah2015}.

Dispersive methods have also played an important role in the study of fluid systems. For small irrotational perturbations of a stable constant state, Audiard and Haspot \cite{AudiardHaspot2017} established global well-posedness for a general class of Euler--Korteweg systems in dimensions \(d\geq3\). Their analysis combines modified energy estimates with dispersive estimates and, in dimensions three and four, a detailed study of the quadratic nonlinear structure through the space--time resonance method. In the particular quantum hydrodynamic regime \(K(\rho)=\kappa/\rho\), global non-vacuum strong solutions were also obtained via the connection with the Gross--Pitaevskii equation \cite{AudiardHaspot2018}, building on the scattering theory for small perturbations of its nonzero equilibrium \cite{GustafsonNakanishiTsai2009}. These results demonstrate the stabilizing role of dispersion in quantum fluid models. The pressureless relativistic system considered here, however, exhibits a different spectral structure, characterized by two interacting dispersive branches and the resulting nontrivial cross-mode resonances.

The stabilizing role of dispersion has also been extensively investigated for the Euler--Poisson and Euler--Maxwell equations. For the three-dimensional electron Euler--Poisson system, Guo \cite{Guo1998} used the Klein--Gordon-type dispersion generated by plasma oscillations, together with a quadratic normal-form transformation, to construct global small smooth irrotational solutions. For the ion dynamics, Guo and Pausader \cite{GuoPausader2011} treated a substantially different dispersion relation by combining refined linear decay estimates, a normal-form transformation, low-frequency analysis, and bilinear Fourier multiplier estimates. Related combinations of dispersive, resonance, and high-order energy methods were subsequently developed for the Euler--Maxwell system and other multi-component plasma models; see, for example, \cite{GermainMasmoudi2014,GuoIonescuPausader2016}. These works illustrate both the stabilizing effect of dispersion and the sensitivity of the nonlinear analysis to the geometry of the underlying dispersion relations and their associated interaction phases.

A distinctive difficulty of the pressureless relativistic quantum hydrodynamic system lies in the simultaneous presence of higher-order spatial and temporal corrections. In addition to the quantum term $n\nabla\left(\frac{\Delta\sqrt n}{\sqrt n}\right)$,
the momentum equation contains the mixed space--time contribution $n\nabla\left(\frac{\partial_{tt}\sqrt n}{\sqrt n}\right)$,
while the continuity equation itself involves the second time derivative of the phase. Thus, the relativistic correction changes the evolution structure of the fluid system rather than entering merely as a lower-order perturbation. To uncover a compatible hyperbolic structure, we introduce the logarithmic amplitude
$\psi=\log\sqrt n$
and reformulate the system as a coupled semilinear wave system. The linear coupling through the time derivatives is skew-symmetric and cancels in the energy identity, leading to a natural high-order energy involving only derivatives of \(\psi\) and \(S\). This derivative structure is particularly well suited to perturbations of a non-vacuum equilibrium, since it yields high-order Sobolev control without requiring \(L^2\)-bounds for the undifferentiated amplitude and phase.

A second structural mechanism governs the global analysis. The linearized system possesses two dispersive branches,
$$
\lambda_\pm(\xi)=\langle\xi\rangle\pm1,
$$
which have the same group velocity but markedly different low-frequency behavior: the slow branch vanishes quadratically at the origin, whereas the fast branch remains separated from zero by a spectral gap. The resulting quadratic mode interactions exhibit a nontrivial resonance geometry. The crucial point is that the nonlinear symbols are exactly compatible with the associated interaction phases. For every active quadratic interaction, the corresponding symbol factorizes as
$$
m(\xi,\eta)=\Phi(\xi,\eta)\widetilde m(\xi,\eta),
$$
where \(\Phi\) denotes the interaction phase and the quotient \(m/\Phi\) remains nonsingular in the relevant frequency regimes. In particular, the quadratic symbol vanishes on the time-resonant set, so that the apparent singularity arising from division by \(\Phi\) in the normal-form transformation is cancelled by the nonlinearity itself. This permits a quadratic normal-form reduction adapted simultaneously to the fast and slow dispersive branches. Combining this cancellation with the derivative energy structure and frequency-localized dispersive estimates, we obtain global control of sufficiently small, smooth, and spatially localized perturbations of the constant non-vacuum equilibrium. Consequently, we establish global existence, persistence of strict positivity of the density, decay of the first derivatives, and scattering in a high-order derivative energy topology.

The remainder of the paper is organized as follows. In Section~\ref{hyperbolic reformulation}, we reformulate the original hydrodynamic system as a coupled hyperbolic system in the logarithmic amplitude and phase variables, revealing its skew-symmetric coupling structure. In Section~\ref{Fourier Analysis}, we analyze the linearized system, identify the two characteristic branches, and establish the relevant dispersive and resonance estimates. In Section~\ref{Quadratic}, we analyze the quadratic nonlinear interactions, establish the phase factorization of the corresponding symbols, and construct the quadratic normal-form transformation. Finally, in Section~\ref{sec:small-data-global}, we combine the nonlinear estimates with the high-order energy estimates and close the bootstrap argument, proving global existence, decay, and scattering.

\section{A hyperbolic reformulation of the RQHD system}\label{hyperbolic reformulation}
In this section, we derive a hyperbolic formulation of the RQHD system and then introduce the logarithmic density variable adapted to the normalized equilibrium.
\subsection{Reduction to a coupled hyperbolic system}
To study the global dynamics of small perturbations of a constant
non-vacuum equilibrium, we first reformulate the relativistic quantum
hydrodynamic system in terms of the logarithmic amplitude and the phase.
We consider the normalized case $\varepsilon=\upsilon=1$.
On any time interval on which the density is strictly positive,
we introduce
\begin{equation}
\omega:=\sqrt{n}>0.
\label{eq:rho-definition}
\end{equation}
We first rewrite the nonlinear transport term as
\begin{align}
\operatorname{div}
\bigl(n\nabla S\otimes\nabla S\bigr)
&=
\frac12 n\nabla|\nabla S|^2
+
\operatorname{div}(n\nabla S)\nabla S
\notag\\
&=
n(\nabla^2S)\nabla S
+
\operatorname{div}(n\nabla S)\nabla S.
\label{eq:momentum-divergence}
\end{align}
Multiplying \(\eqref{J1}_1\) by \(\nabla S\) and substituting the
resulting identity into \(\eqref{J1}_2\), we obtain
\begin{equation}
\partial_t S
+\frac12|\nabla S|^2
-\frac12\frac{\Delta\omega}{\omega}
=
\frac12(\partial_tS)^2
-\frac12\frac{\partial_{tt}\omega}{\omega}.
\label{eq:phase-omega-form}
\end{equation}
Equivalently,
\begin{equation}
\partial_{tt}\omega-\Delta\omega
=
\omega
\left(
(\partial_tS)^2
-2\partial_tS
-|\nabla S|^2
\right).
\label{eq:omega-wave}
\end{equation}
On the other hand, rewriting \(\eqref{J1}_1\) in terms of
\(\omega=\sqrt n\) gives
\begin{equation}
\partial_{tt}S-\Delta S
=
\frac{2}{\omega}
\left(
\partial_t\omega
+\nabla\omega\cdot\nabla S
-\partial_t\omega\,\partial_tS
\right).
\label{eq:S-wave}
\end{equation}
Hence, as long as \(\omega>0\), the RQHD system can be reformulated as
the coupled hyperbolic system
\begin{equation}
\left\{
\begin{aligned}
\partial_{tt}\omega-\Delta\omega
&=
\omega
\left(
(\partial_tS)^2
-2\partial_tS
-|\nabla S|^2
\right),
\\[1mm]
\partial_{tt}S-\Delta S
&=
\frac{2}{\omega}
\left(
\partial_t\omega
+\nabla\omega\cdot\nabla S
-\partial_t\omega\,\partial_tS
\right).
\end{aligned}
\right.
\label{eq:omega-S-hyperbolic-system}
\end{equation}
The corresponding initial data are
\begin{equation}
(\omega,\partial_t\omega)(0,x)
=
(\omega_0,\omega_1)(x)
:=
\left(
\sqrt{n_0(x)},
\frac{n_1(x)}{2\sqrt{n_0(x)}}
\right),
\qquad x\in\mathbb R^3,
\label{eq:omega-initial-data}
\end{equation}
and
\begin{equation}
(S,\partial_tS)(0,x)
=
(S_0,S_1)(x),
\qquad x\in\mathbb R^3,
\label{eq:S-initial-data}
\end{equation}
where
\[
\inf_{x\in\mathbb R^3}\omega_0(x)>0.
\]

\subsection{Logarithmic density and the normalized equilibrium}
We now introduce the logarithmic amplitude variable
\begin{equation}
\psi:=\log\omega=\frac12\log n.
\label{eq:log-amplitude}
\end{equation}
Since \(\omega=e^\psi\), we have
\begin{align}
\partial_t\omega
&=\omega\,\partial_t\psi,
&
\nabla\omega
&=\omega\,\nabla\psi,
\\
\partial_{tt}\omega
&=
\omega
\left(
\partial_{tt}\psi+(\partial_t\psi)^2
\right),
&
\Delta\omega
&=
\omega
\left(
\Delta\psi+|\nabla\psi|^2
\right).
\end{align}
Substituting these identities into
\eqref{eq:omega-S-hyperbolic-system}, we obtain
\begin{equation}
\left\{
\begin{aligned}
\partial_{tt}\psi-\Delta\psi+2\partial_tS
&=
(\partial_tS)^2
-|\nabla S|^2
-(\partial_t\psi)^2
+|\nabla\psi|^2,
\\[1mm]
\partial_{tt}S-\Delta S-2\partial_t\psi
&=
2\nabla\psi\cdot\nabla S
-2\partial_t\psi\,\partial_tS.
\end{aligned}
\right.
\label{eq:psi-S-system}
\end{equation}
Thus the RQHD system is reduced to a coupled semilinear wave system for
\((\psi,S)\), whose nonlinearities are quadratic in first-order
derivatives.
The corresponding initial data are
\begin{equation}
(\psi,\partial_t\psi)(0,x)
=
(\psi_0,\psi_1)(x)
:=
\left(
\frac12\log n_0(x),
\frac{n_1(x)}{2n_0(x)}
\right),
\label{eq:psi-initial-data}
\end{equation}
and
\begin{equation}
(S,\partial_tS)(0,x)
=
(S_0,S_1)(x),
\qquad x\in\mathbb R^3,
\label{eq:S-initial-data}
\end{equation}
with
\[
\inf_{x\in\mathbb R^3} n_0(x)>0.
\]
The hyperbolic formulation above provides the natural framework for the local theory. Our main interest, however, is the long-time behavior of solutions near the constant non-vacuum equilibrium. This leads to the following small-data global existence problem.
\begin{problem}[Small-data global existence]\label{prob:main}
Consider the Cauchy problem for the coupled wave-type system
\begin{align}\label{eq:main-system}
\begin{cases}
\psi_{tt}-\Delta\psi
+2S_t
=
S_t^2-|\nabla S|^2-\psi_t^2+|\nabla\psi|^2,\\[0.6ex]
S_{tt}-\Delta S
-2\psi_t
=
2\nabla\psi\cdot\nabla S
-2\psi_tS_t,
\end{cases}
\end{align}
in \((0,\infty)\times\mathbb R^3\), with initial data
\begin{equation}
\begin{aligned}
&\psi(0,x)=\delta\psi_0(x),\qquad
\partial_t\psi(0,x)=\delta\psi_1(x),\\
&S(0,x)=\delta S_0(x),\qquad
\partial_tS(0,x)=\delta S_1(x).
\end{aligned}
\label{eq:main-data}
\end{equation}
Here \(\delta>0\) is a small amplitude parameter.
For fixed initial profiles, the problem is to determine whether
there exists \(\delta_*>0\) such that, for every
\(0<\delta\leq\delta_*\), the Cauchy problem
\eqref{eq:main-system}--\eqref{eq:main-data}
admits a unique global classical solution.
\end{problem}
\begin{remark}
The logarithmic transformation is well defined as long as the density remains strictly positive. It eliminates the reciprocal density factors and reveals the quadratic derivative structure of the system.
\end{remark}
\paragraph{Proof strategy.}
We first diagonalize the linearized system into four half-wave modes with frequencies
$$
\pm\lambda_\pm(\xi),
\qquad
\lambda_\pm(\xi)=\langle\xi\rangle\pm1.
$$
The main low-frequency difficulty comes from the slow branch
$$
\lambda_-(\xi)\sim \frac12|\xi|^2
\qquad (|\xi|\ll1).
$$
We therefore work with modal derivative norms adapted to the two branches. These norms are equivalent to the natural derivative energy of \((\psi,S)\) and avoid the need for separate \(L^2\)-control of the undifferentiated variables.
For the nonlinear problem, every active quadratic interaction satisfies the factorization
$$
m=\Phi\,\widetilde m,
$$
where \(\Phi\) is the interaction phase and \(\widetilde m\) is regular. Hence a quadratic normal-form transformation removes all quadratic terms without introducing a singular divisor, and the transformed equation contains only cubic nonlinearities. Dyadic dispersive estimates, together with the spatially weighted localization of the initial data, yield the required \(\langle t\rangle^{-3/2}\) decay. This integrable decay is then inserted into the high-order energy estimate to close the global bootstrap. Finally, the time integrability of the cubic remainder and the decay of the quadratic normal-form correction give scattering in the derivative energy norm.

\section{Fourier Analysis of the Linearized System}\label{Fourier Analysis}
In this section, we analyze the linearized system through its Fourier spectral decomposition and compare the modal derivative norms with the natural derivative norms of the original variables, with particular attention to the low-frequency degeneracy of the slow branch. We then establish the corresponding dispersive decay estimates. Finally, we analyze the interaction phases and characterize the associated time, space, and space-time resonant sets.
\subsection{Spectral Structure of the Linearized System}
We consider the linearized system
\begin{equation}
\label{eq:linearized-system}
\begin{cases}
\widetilde{\psi}_{tt}-\Delta\widetilde{\psi}
+2\widetilde{S}_t=0,\\[0.4em]
\widetilde{S}_{tt}-\Delta\widetilde{S}
-2\widetilde{\psi}_t=0,
\end{cases}
\qquad
(t,x)\in[0,\infty)\times\mathbb R^3.
\end{equation}
To determine its characteristic frequencies, we take the Fourier
transform in the spatial variable \(x\), using the convention
\[
\widehat f(t,\xi)
=
\mathcal F_x[f](t,\xi)
=
\int_{\mathbb R^3}e^{-ix\cdot\xi}f(t,x)\,dx.
\]
For simplicity, throughout this linear analysis we write
\[
\widehat\psi:=\mathcal F_x[\widetilde{\psi}],
\qquad
\widehat S:=\mathcal F_x[\widetilde{S}].
\]
Since the spatial Fourier transform commutes with time
differentiation and satisfies
\(\mathcal F_x[-\Delta f]=|\xi|^2\widehat f\),
system \eqref{eq:linearized-system} becomes
\begin{equation}
\label{eq:linearized-fourier-system}
\begin{cases}
\partial_t^2\widehat\psi
+|\xi|^2\widehat\psi
+2\partial_t\widehat S=0,\\[0.4em]
\partial_t^2\widehat S
+|\xi|^2\widehat S
-2\partial_t\widehat\psi=0.
\end{cases}
\end{equation}
For each fixed \(\xi\in\mathbb R^3\), this is a coupled system
of ordinary differential equations in time.
We next diagonalize it to identify the characteristic frequencies
and the corresponding linear modes.
Introducing
\[
Y(t,\xi)
:=
\begin{pmatrix}
\widehat\psi(t,\xi)\\
\widehat S(t,\xi)\\
\partial_t\widehat\psi(t,\xi)\\
\partial_t\widehat S(t,\xi)
\end{pmatrix},
\]
we rewrite \eqref{eq:linearized-fourier-system} as
\begin{equation}\label{eq:first-order-fourier-system}
\partial_tY(t,\xi)=A(\xi)Y(t,\xi),
\end{equation}
where
\[
A(\xi)
=
\begin{pmatrix}
0&0&1&0\\
0&0&0&1\\
-|\xi|^2&0&0&-2\\
0&-|\xi|^2&2&0
\end{pmatrix}.
\]
Let
\[
\Lambda(\xi):=\sqrt{1+|\xi|^2},
\qquad
\lambda_\sigma(\xi):=\Lambda(\xi)+\sigma,
\qquad
\sigma\in\{\pm1\}.
\]
A direct computation gives
\[
\det\bigl(\tau I-A(\xi)\bigr)
=
(\tau^2+|\xi|^2)^2+4\tau^2
=
\prod_{\sigma\in\{\pm1\}}
\left(\tau^2+\lambda_\sigma(\xi)^2\right).
\]
Hence
\[
\operatorname{spec}A(\xi)
=
\left\{
\varepsilon i\lambda_\sigma(\xi):
\sigma,\varepsilon\in\{\pm1\}
\right\}.
\]
Here $\sigma$ labels the two dispersion branches, while $\varepsilon$
gives the two time directions. At low frequency,
\[
\lambda_-(\xi)
=
\frac12|\xi|^2+O(|\xi|^4),
\qquad
\lambda_+(\xi)
=
2+\frac12|\xi|^2+O(|\xi|^4).
\]
Thus the lower branch is Schr\"odinger-like near $\xi=0$, whereas the
upper branch has a nonzero frequency gap. At high frequencies,
\[
\lambda_\pm(\xi)
=
|\xi|\pm1+O(|\xi|^{-1}),
\]
so both branches are wave-like. Although the shifts $\pm1$
do not affect the linear group velocities, they enter the nonlinear
interaction phases and therefore affect the time-resonance structure.
For each $\sigma,\varepsilon\in\{\pm1\}$, an associated eigenvector is
\[
v_\sigma^\varepsilon(\xi)
=
\begin{pmatrix}
1\\
-\varepsilon\sigma i\\
\varepsilon i\lambda_\sigma(\xi)\\
\sigma\lambda_\sigma(\xi)
\end{pmatrix},
\]
and
\[
A(\xi)v_\sigma^\varepsilon(\xi)
=
\varepsilon i\lambda_\sigma(\xi)
v_\sigma^\varepsilon(\xi).
\]
Consequently,
\[
Y(t,\xi)
=
\sum_{\sigma,\varepsilon\in\{\pm1\}}
e^{\varepsilon it\lambda_\sigma(\xi)}
a_{\sigma,0}^\varepsilon(\xi)
v_\sigma^\varepsilon(\xi),
\]
where $a_{\sigma,0}^\varepsilon$ are determined by the initial data.

\begin{remark}
Although system \eqref{eq:main-system} is equivalent to a Klein--Gordon equation after a nonlinear change of variables and removal of the rest-mass oscillation, its linearized spectrum in the variables $(\psi,S)$ consists of two shifted branches. In particular, the slow branch is gapless and behaves quadratically near the origin, while both branches are wave-like at high frequencies.
\end{remark}
\subsection{Equivalence of Derivative Norms}
We use the eigenvector basis above to define modal coordinates
for the nonlinear solution and relate their derivative norms
to those of $\psi$ and $S$.
For $\alpha=(\sigma,\varepsilon)$, the inverse modal transformation
gives
\begin{equation}
\widehat U_\sigma^\varepsilon
=
\frac{1}{4\Lambda}
\left[
\lambda_{-\sigma}
\bigl(\widehat\psi+i\varepsilon\sigma\widehat S\bigr)
+\sigma\widehat{S_t}-i\varepsilon\widehat{\psi_t}
\right].
\label{eq:modal-coordinate-definition}
\end{equation}
The reconstruction formulas are
\begin{align}
\widehat\psi
&=\sum_{\sigma,\varepsilon}\widehat U_\sigma^\varepsilon,
&
\widehat S
&=-i\sum_{\sigma,\varepsilon}
\varepsilon\sigma\widehat U_\sigma^\varepsilon,
\label{eq:modal-reconstruction-position}
\\
\widehat{\psi_t}
&=i\sum_{\sigma,\varepsilon}
\varepsilon\lambda_\sigma\widehat U_\sigma^\varepsilon,
&
\widehat{S_t}
&=\sum_{\sigma,\varepsilon}
\sigma\lambda_\sigma\widehat U_\sigma^\varepsilon.
\label{eq:modal-reconstruction-velocity}
\end{align}
All sums over $\alpha$ run over these four pairs.
Define
\begin{equation}
\mathcal D_\alpha h
:=
\bigl(\nabla h,\lambda_\sigma(D)h\bigr),
\qquad
\|\mathcal D_\alpha h\|_X
:=
\|\nabla h\|_X+\|\lambda_\sigma(D)h\|_X.
\label{eq:first-order-modal-operator}
\end{equation}
The corresponding modal space is
\begin{equation}
\mathcal H_\sigma^s
:=
\{h:\nabla h,\lambda_\sigma(D)h\in H^s\},
\qquad
\|h\|_{\mathcal H_\sigma^s}
:=\|\mathcal D_\alpha h\|_{H^s}.
\label{eq:homogeneous-modal-space}
\end{equation}
For $\sigma=-1$, the space is understood modulo additive constants.
Whenever undifferentiated slow modes occur in products, we use their
unique $L^6(\mathbb R^3)$ representatives, supplied by the homogeneous
Sobolev embedding. For $\sigma=+1$, the bound $\lambda_+\geq2$ already
implies $h\in H^s$, so no quotient by constants is needed.

\begin{lemma}[Equivalence of derivative norms]
\label{lem:modal-energy-equivalence}
For every real $s\geq0$, the modal coordinates defined by
\eqref{eq:modal-coordinate-definition} satisfy
\begin{equation}
\begin{aligned}
\sum_\alpha
\|\mathcal D_\alpha U_\alpha\|_{H^s}
\simeq{}&
\|\nabla\psi\|_{H^s}+\|\nabla S\|_{H^s}+\|\psi_t\|_{H^s}+\|S_t\|_{H^s}.
\end{aligned}
\label{eq:modal-physical-equivalence}
\end{equation}
More generally, for every nonnegative integer $L$,
\begin{equation}
\begin{aligned}
&\sum_\alpha\sum_{|\rho|\leq L}
\left(
\|x^\rho\nabla U_\alpha\|_{H^s}
+\|x^\rho\lambda_\sigma(D)U_\alpha\|_{H^s}
\right)
\\
&\qquad\simeq
\sum_{|\rho|\leq L}
\Bigl(
\|x^\rho\nabla\psi\|_{H^s}
+\|x^\rho\nabla S\|_{H^s}
+\|x^\rho\psi_t\|_{H^s}
+\|x^\rho S_t\|_{H^s}
\Bigr).
\end{aligned}
\label{eq:weighted-modal-physical-equivalence}
\end{equation}
The constants may depend on $s$ and $L$. The weights in
\eqref{eq:weighted-modal-physical-equivalence} are applied after
the derivatives.
\end{lemma}

\begin{proof}
The reconstruction formulas
\eqref{eq:modal-reconstruction-position}--%
\eqref{eq:modal-reconstruction-velocity} immediately give
\[
\|\nabla\psi\|_{H^s}+\|\nabla S\|_{H^s}
+\|\psi_t\|_{H^s}+\|S_t\|_{H^s}
\lesssim
\sum_\alpha
\|\mathcal D_\alpha U_\alpha\|_{H^s}.
\]
The same argument applies after multiplying each differentiated
quantity by $x^\rho$.
To prove the reverse inequalities, we express the modal derivatives
using only derivatives of $\psi$ and $S$. Set
\[
a_\sigma(\xi):=\frac{\lambda_{-\sigma}(\xi)}{\Lambda(\xi)}
=1-\frac{\sigma}{\Lambda(\xi)},
\qquad
c_\sigma(\xi):=\frac{\lambda_\sigma(\xi)}{\Lambda(\xi)}
=1+\frac{\sigma}{\Lambda(\xi)},
\qquad
T_j:=\langle D\rangle^{-1}\partial_j.
\]
Differentiating \eqref{eq:modal-coordinate-definition} gives
\begin{equation}
\begin{aligned}
4\partial_jU_\sigma^\varepsilon
={}&
a_\sigma(D)\partial_j\psi
+i\varepsilon\sigma a_\sigma(D)\partial_jS+\sigma T_jS_t-i\varepsilon T_j\psi_t.
\end{aligned}
\label{eq:modal-gradient-from-derivatives}
\end{equation}
Moreover, the identity
\begin{equation}
\lambda_+(\xi)\lambda_-(\xi)=|\xi|^2
\label{eq:lambda-product-identity}
\end{equation}
and the relation
$|D|^2\langle D\rangle^{-1}=-\sum_{j=1}^3T_j\partial_j$
yield
\begin{equation}
\begin{aligned}
4\lambda_\sigma(D)U_\sigma^\varepsilon
={}&
-\sum_{j=1}^3T_j
\bigl(\partial_j\psi+i\varepsilon\sigma\partial_jS\bigr)
+\sigma c_\sigma(D)S_t
-i\varepsilon c_\sigma(D)\psi_t.
\end{aligned}
\label{eq:modal-lambda-from-derivatives}
\end{equation}
The symbols $a_\sigma$, $c_\sigma$, and $i\xi_j/\Lambda(\xi)$
are smooth symbols of order zero. In particular, their multipliers
are bounded on $H^s$. Taking norms in
\eqref{eq:modal-gradient-from-derivatives} and
\eqref{eq:modal-lambda-from-derivatives} proves the reverse inequality
in \eqref{eq:modal-physical-equivalence}.
For the weighted estimate, let $a$ be any of these order-zero symbols.
For smooth functions, the Fourier product rule gives
\begin{equation}
x^\rho a(D)g
=
\sum_{\theta\leq\rho}
\binom{\rho}{\theta}i^{|\theta|}
(\partial_\xi^\theta a)(D)
\bigl(x^{\rho-\theta}g\bigr).
\label{eq:weighted-order-zero-commutator}
\end{equation}
Every differentiated symbol on the right is bounded and therefore
defines a bounded multiplier on $H^s$. Hence
\[
\sum_{|\rho|\leq L}\|x^\rho a(D)g\|_{H^s}
\lesssim_{s,L}
\sum_{|\rho|\leq L}\|x^\rho g\|_{H^s}.
\]
Apply this estimate to
\eqref{eq:modal-gradient-from-derivatives} and
\eqref{eq:modal-lambda-from-derivatives}, where $g$ is a component
of $\nabla\psi$, $\nabla S$, $\psi_t$, or $S_t$.
Together with the weighted reconstruction bound, this proves
\eqref{eq:weighted-modal-physical-equivalence}.
The estimates extend to the stated spaces by approximation.
In particular, no separate $L^2$ assumption on $\psi$ or $S$
is used.
\end{proof}

For the linearized system, the modal equations are
\begin{equation}
(\partial_t-i\varepsilon\lambda_\sigma(D))
U_\sigma^\varepsilon=0.
\label{eq:modal-linear-equation}
\end{equation}
For either the linear or nonlinear evolution, define the profiles by
\begin{equation}
f_\alpha(t)
:=
e^{-it\omega_\alpha(D)}U_\alpha(t),
\qquad
\omega_\alpha(\xi):=\varepsilon\lambda_\sigma(\xi).
\label{eq:profile-definition}
\end{equation}
The profiles are constant in time for the linearized system.
Lemma~\ref{lem:modal-energy-equivalence} allows estimates for the
modal derivatives to be transferred to the derivatives of $\psi$
and $S$ used in the energy and scattering statements.
\begin{remark}
Lemma~\ref{lem:modal-energy-equivalence} connects the modal
estimates with the energy of the original system.
It shows that controlling
$\nabla U_\alpha$ and $\lambda_\sigma(D)U_\alpha$
is equivalent to controlling
$\nabla\psi$, $\nabla S$, $\psi_t$, and $S_t$.
Although the slow branch satisfies $\lambda_-(0)=0$,
this equivalence remains uniform at low frequencies.
Thus the modal analysis requires no additional
$L^2$ assumption on the undifferentiated variables
$\psi$ and $S$.
\end{remark}
\subsection{Dyadic Dispersive Estimates}
We now estimate the linear evolution of each modal component.
Recall that
\[
\lambda_\sigma(\xi)=\sqrt{1+|\xi|^2}+\sigma,
\qquad \sigma\in\{\pm1\},
\]
where $\xi\in\mathbb R^3$ is the spatial Fourier variable and
$\sigma$ labels the dispersion branch. The sign
$\varepsilon\in\{\pm1\}$ selects the temporal frequency
$\varepsilon\lambda_\sigma(\xi)$. The corresponding linear
propagator is defined by
\[
\mathcal F_x\!\left[e^{\varepsilon it\lambda_\sigma(D)}f\right](\xi)
=e^{\varepsilon it\lambda_\sigma(\xi)}\widehat f(\xi),
\qquad t\in\mathbb R.
\]
Let $(P_k)_{k\in\mathbb Z}$ be the homogeneous Littlewood--Paley
projections, with
\[
\widehat{P_kf}(\xi)=\varphi(2^{-k}\xi)\widehat f(\xi),
\]
where $\varphi\in C_c^\infty(\mathbb R^3)$ is a fixed smooth radial
cutoff supported in $\{\xi:1/2\leq|\xi|\leq2\}$.
Thus $k$ indexes the frequency blocks, and $P_kf$ contains
frequencies comparable to $2^k$.
\begin{proposition}[Dyadic dispersive estimate]
\label{prop:dyadic-dispersive}
For every $f\in L^1(\mathbb R^3)$, $t\in\mathbb R$,
$k\in\mathbb Z$, and $\sigma,\varepsilon\in\{\pm1\}$,
\[
\left\|
e^{\varepsilon it\lambda_\sigma(D)}P_kf
\right\|_{L^\infty(\mathbb R^3)}
\lesssim a_k(t)\|f\|_{L^1(\mathbb R^3)},
\]
where the scalar decay factor $a_k(t)$ is defined by
\[
a_k(t):=
\begin{cases}
2^{3k}(1+|t|2^{2k})^{-3/2},
& k\leq0,\\[0.5em]
\displaystyle
\min\left\{
2^{3k},\ |t|^{-1}2^{2k},\ |t|^{-3/2}2^{5k/2}
\right\},
& k\geq1,\ t\neq0.
\end{cases}
\]
For $k\geq1$, we set $a_k(0)=2^{3k}$.
The implicit constant depends only on the fixed cutoff $\varphi$
and is independent of $k$, $t$, $\sigma$, and $\varepsilon$.
In addition, let $\chi\in C_c^\infty(\mathbb R^3)$ be a fixed
smooth compactly supported frequency cutoff, and define
$P_{\mathrm{comp}}=\chi(D)$ by
\[
\widehat{P_{\mathrm{comp}}f}(\xi)
=\chi(\xi)\widehat f(\xi).
\]
Then
\[
\left\|
e^{\varepsilon it\lambda_\sigma(D)}P_{\mathrm{comp}}f
\right\|_{L^\infty(\mathbb R^3)}
\leq C_\chi(1+|t|)^{-3/2}\|f\|_{L^1(\mathbb R^3)},
\]
where $C_\chi$ depends on $\chi$ but is independent of $t$,
$\sigma$, and $\varepsilon$.
\end{proposition}
\begin{proof}
Since $\lambda_\sigma(\xi)=\langle\xi\rangle+\sigma$,
\[
e^{\varepsilon it\lambda_\sigma(D)}
=e^{\varepsilon i\sigma t}
e^{\varepsilon it\langle D\rangle}.
\]
The scalar factor $e^{\varepsilon i\sigma t}$ has modulus one,
so it suffices to estimate the Klein--Gordon propagator
$e^{\varepsilon it\langle D\rangle}$.
For its radial phase $\phi(r)=\sqrt{1+r^2}$, where $r=|\xi|$,
\[
\phi'(r)=\frac{r}{\sqrt{1+r^2}},
\qquad
\phi''(r)=(1+r^2)^{-3/2}.
\]
At low frequencies, $\phi'(r)\simeq r$ and
$\phi''(r)\simeq1$; at high frequencies,
$\phi'(r)\simeq1$ and $\phi''(r)\simeq r^{-3}$.
The dyadic estimates therefore follow from
\cite[Theorem~1(a)--(b) and Section~4]{GuoPengWang2008},
together with the elementary kernel bound $2^{3k}$.
For the compact-frequency estimate, we use the low-frequency
bound in \cite[Theorem~1(c)]{GuoPengWang2008} and the dyadic
estimates for the finitely many higher-frequency blocks meeting
the support of $\chi$. The smooth cutoff $\chi(D)$ is bounded
on $L^1$, with a bound depending only on $\chi$.
This proves the stated estimate for $P_{\mathrm{comp}}$.
\end{proof}

\subsection{Quadratic Space--Time Resonances}

For an output mode $(\sigma,\varepsilon)$ and input modes
$(\mu,\varepsilon_1)$ and $(\nu,\varepsilon_2)$, the interaction
phase in the profile equation is
\[
-\varepsilon\lambda_\sigma(\xi)
+\varepsilon_1\lambda_\mu(\eta)
+\varepsilon_2\lambda_\nu(\xi-\eta).
\]
Since $\lambda_\rho(\zeta)=\langle\zeta\rangle+\rho$, its constant
shift is
\begin{equation}
h=-\varepsilon\sigma+\varepsilon_1\mu+\varepsilon_2\nu
\in\{-3,-1,1,3\}.
\label{eq:general-phase-shift}
\end{equation}
We first classify the space--time resonances of all these phases.
The quadratic symbols computed in the next section determine which
interactions actually occur.

\begin{proposition}[Classification of quadratic space--time resonances]
\label{prop:space-time-resonance}
For $h\in\{-3,-1,1,3\}$ and
$\varepsilon,\varepsilon_1,\varepsilon_2\in\{\pm1\}$, let
\[
\Phi_{h;\varepsilon,\varepsilon_1,\varepsilon_2}(\xi,\eta)
=h-\varepsilon\langle\xi\rangle
+\varepsilon_1\langle\eta\rangle
+\varepsilon_2\langle\xi-\eta\rangle.
\]
Define
\[
\mathcal T=\{(\xi,\eta):\Phi(\xi,\eta)=0\},
\qquad
\mathcal S=\{(\xi,\eta):\nabla_\eta\Phi(\xi,\eta)=0\},
\qquad
\mathcal R=\mathcal T\cap\mathcal S.
\]
The space-resonant set is
\[
\mathcal S_{h;\varepsilon,\varepsilon_1,\varepsilon_2}
=
\begin{cases}
\{(\xi,\eta):\eta=\xi/2\},
&\varepsilon_1=\varepsilon_2,\\[1mm]
\{(0,\eta):\eta\in\mathbb R^3\},
&\varepsilon_1=-\varepsilon_2.
\end{cases}
\]
The only nonempty space--time resonant sets are
\begin{align*}
\mathcal R_{h;h,a,-a}
&=\{(0,\eta):\eta\in\mathbb R^3\},
&&h,a\in\{\pm1\},\\
\mathcal R_{h;-h,-h,-h}
&=\{(0,0)\},
&&h\in\{\pm1\},\\
\mathcal R_{3a;a,-a,-a}
&=\{(0,0)\},
&&a\in\{\pm1\}.
\end{align*}
All other space--time resonant sets are empty.
\end{proposition}

\begin{proof}
Differentiating the phase gives
\[
\nabla_\eta\Phi
=\varepsilon_1\frac{\eta}{\langle\eta\rangle}
-\varepsilon_2\frac{\xi-\eta}{\langle\xi-\eta\rangle}.
\]
The map $v\mapsto v/\sqrt{1+|v|^2}$ is odd and injective.
Thus $\nabla_\eta\Phi=0$ implies $\eta=\xi-\eta$ when
$\varepsilon_1=\varepsilon_2$, and $\eta=-(\xi-\eta)$ when
$\varepsilon_1=-\varepsilon_2$. The converse implications are
immediate, proving the description of $\mathcal S$.
If $\varepsilon_1=-\varepsilon_2$, then $\xi=0$ on
$\mathcal S$ and
\[
\Phi(0,\eta)=h-\varepsilon.
\]
Hence $\mathcal R=\{(0,\eta):\eta\in\mathbb R^3\}$ precisely
when $h=\varepsilon\in\{\pm1\}$, and is empty otherwise.
If $\varepsilon_1=\varepsilon_2=a$, then $\eta=\xi/2$.
Set
\[
r=|\xi|,
\qquad C_r=\sqrt{1+r^2},
\qquad D_r=\sqrt{4+r^2}.
\]
The time-resonance condition becomes
\begin{equation}
h=\varepsilon C_r-aD_r.
\label{eq:resonance-radial-equation}
\end{equation}
When $\varepsilon=a$, we have
\[
h=-\varepsilon(D_r-C_r),
\qquad
D_r-C_r=\frac{3}{D_r+C_r}\in(0,1].
\]
Since $h\in\{-3,-1,1,3\}$, equality is possible only when
$r=0$ and $h=-\varepsilon=-a$. This gives
$\mathcal R_{h;-h,-h,-h}=\{(0,0)\}$ for $h=\pm1$.
When $\varepsilon=-a$, we instead have
\[
h=\varepsilon(C_r+D_r),
\qquad C_r+D_r\geq3,
\]
with equality only at $r=0$. Thus the only possibility is
$h=3\varepsilon$ and $r=0$, giving
$\mathcal R_{3\varepsilon;\varepsilon,-\varepsilon,-\varepsilon}
=\{(0,0)\}$.
\end{proof}
The cases $h=\pm3$ are included here as possible phase resonances.
We will show below that their quadratic symbols vanish identically,
so they do not contribute to the nonlinear evolution.
\section{Quadratic Fourier structure}\label{Quadratic}
In this section, we analyze the quadratic modal interactions and establish the exact cancellation of the time-resonant phases. This structure allows us to perform a nonsingular quadratic normal-form reduction, for which we then derive the required dyadic multiplier estimates.
\subsection{Quadratic Fourier Symbols and Modal Interactions}
The quadratic nonlinearities are
\begin{align}
\mathcal N_1(\psi,S)
&=S_t^2-|\nabla S|^2-\psi_t^2+|\nabla\psi|^2,
\label{eq:quadratic-N1}
\\
\mathcal N_2(\psi,S)
&=2\nabla\psi\cdot\nabla S-2\psi_tS_t.
\label{eq:quadratic-N2}
\end{align}
We use the convolution normalization
\[
\widehat{uv}(\xi)
=\frac{1}{(2\pi)^3}\int_{\mathbb R^3}
\widehat u(\eta)\widehat v(\xi-\eta)\,d\eta.
\]
Spatial derivatives contribute Fourier frequencies, while the
dispersion factors in the time-derivative terms come from the modal
reconstruction of $\psi_t$ and $S_t$.
For input modes $(\mu,\varepsilon_1)$ and $(\nu,\varepsilon_2)$, set
\[
\zeta=\xi-\eta,
\qquad
\lambda_\mu(\eta)=\langle\eta\rangle+\mu,
\qquad
\lambda_\nu(\zeta)=\langle\zeta\rangle+\nu,
\]
and define
\begin{equation}
q_{\mu\nu}(\eta,\zeta)
:=\mu\nu\lambda_\mu(\eta)\lambda_\nu(\zeta)
-\eta\cdot\zeta.
\label{eq:quadratic-null-symbol}
\end{equation}
Substituting the modal reconstruction formulas gives
\begin{align}
\widehat{\mathcal N_1}(\xi)
={}&\frac{1}{(2\pi)^3}
\sum_{\mu,\nu}\sum_{\varepsilon_1,\varepsilon_2}
\int_{\mathbb R^3}
\bigl(1+\varepsilon_1\varepsilon_2\mu\nu\bigr)
q_{\mu\nu}(\eta,\xi-\eta)
\notag\\
&\qquad\times
\widehat U_\mu^{\varepsilon_1}(\eta)
\widehat U_\nu^{\varepsilon_2}(\xi-\eta)\,d\eta,
\label{eq:N1-modal-fourier}
\\
\widehat{\mathcal N_2}(\xi)
={}&-\frac{i}{(2\pi)^3}
\sum_{\mu,\nu}\sum_{\varepsilon_1,\varepsilon_2}
\int_{\mathbb R^3}
\bigl(\varepsilon_1\mu+\varepsilon_2\nu\bigr)
q_{\mu\nu}(\eta,\xi-\eta)
\notag\\
&\qquad\times
\widehat U_\mu^{\varepsilon_1}(\eta)
\widehat U_\nu^{\varepsilon_2}(\xi-\eta)\,d\eta.
\label{eq:N2-modal-fourier}
\end{align}
Here \eqref{eq:N2-modal-fourier} is written in symmetrized form,
obtained by exchanging the two input modes and their frequencies
and averaging. All discrete indices range over $\{\pm1\}$,
and the time variable is suppressed.
The coefficients of $\widehat{\psi_t}$ and $\widehat{S_t}$ in the
inverse modal map are $-i\varepsilon/(4\langle\xi\rangle)$ and
$\sigma/(4\langle\xi\rangle)$, respectively.
Projecting the source onto the output mode $(\sigma,\varepsilon)$
therefore gives
\begin{equation}
\left(\partial_t-i\varepsilon\lambda_\sigma(D)\right)
U_\sigma^\varepsilon
=\frac{1}{4\langle D\rangle}
\left(\sigma\mathcal N_2-i\varepsilon\mathcal N_1\right).
\label{eq:modal-forced-equation}
\end{equation}
Thus the full bilinear symbol is
\begin{equation}
\begin{aligned}
m_{\sigma;\mu,\nu}^{\varepsilon;\varepsilon_1,\varepsilon_2}
(\xi,\eta)
={}&-\frac{i}{4\langle\xi\rangle}
\Bigl[\sigma(\varepsilon_1\mu+\varepsilon_2\nu)
\\
&\qquad+\varepsilon
\bigl(1+\varepsilon_1\varepsilon_2\mu\nu\bigr)\Bigr]
q_{\mu\nu}(\eta,\xi-\eta),
\end{aligned}
\label{eq:full-bilinear-symbol}
\end{equation}
and the modal equation takes the form
\begin{align}
&\left(\partial_t-i\varepsilon\lambda_\sigma(\xi)\right)
\widehat U_\sigma^\varepsilon(\xi)
\notag\\
&\quad=\frac{1}{(2\pi)^3}
\sum_{\mu,\nu}\sum_{\varepsilon_1,\varepsilon_2}
\int_{\mathbb R^3}
m_{\sigma;\mu,\nu}^{\varepsilon;\varepsilon_1,\varepsilon_2}
(\xi,\eta)
\widehat U_\mu^{\varepsilon_1}(\eta)
\widehat U_\nu^{\varepsilon_2}(\xi-\eta)\,d\eta.
\label{eq:full-quadratic-modal-interaction}
\end{align}
The common factor $q_{\mu\nu}$ encodes the differential structure
of the quadratic sources. To determine which modal interactions
are present, set
\[
r_0=\varepsilon\sigma,
\qquad
r_1=\varepsilon_1\mu,
\qquad
r_2=\varepsilon_2\nu.
\]
The discrete coefficient in \eqref{eq:full-bilinear-symbol}
factors as
\[
\sigma(r_1+r_2)+\varepsilon(1+r_1r_2)
=
\varepsilon(1+r_0r_1)(1+r_0r_2).
\]
It is nonzero precisely when $r_0=r_1=r_2$, in which case
its value is $4\varepsilon$. Thus the active interactions satisfy
\[
\varepsilon\sigma
=
\varepsilon_1\mu
=
\varepsilon_2\nu.
\]
For such interactions, the constant phase shift becomes
\[
h=-r_0+r_1+r_2=r_0\in\{\pm1\}.
\]
In particular, all interactions with $h=\pm3$ have identically
vanishing quadratic symbols. Their space--time resonances,
although included in the preceding classification, therefore
do not contribute to the quadratic nonlinearity.

\subsection{Exact Cancellation of the Quadratic Time Resonances}
The following proposition identifies the active interactions and
shows that each corresponding quadratic symbol contains the
interaction phase as an exact factor.
\begin{proposition}[Exact phase factorization]
\label{prop:exact-phase-factorization}
Let
\[
\Phi_{\sigma;\mu,\nu}^{\varepsilon;\varepsilon_1,\varepsilon_2}
(\xi,\eta)
=-\varepsilon\lambda_\sigma(\xi)
+\varepsilon_1\lambda_\mu(\eta)
+\varepsilon_2\lambda_\nu(\xi-\eta).
\]
The quadratic symbol vanishes identically unless there exists
$h\in\{\pm1\}$ such that
\begin{equation}
\varepsilon\sigma=\varepsilon_1\mu=\varepsilon_2\nu=h.
\label{eq:active-mode-condition}
\end{equation}
Equivalently,
\[
\sigma=h\varepsilon,
\qquad\mu=h\varepsilon_1,
\qquad\nu=h\varepsilon_2.
\]
For every active interaction satisfying
\eqref{eq:active-mode-condition}, one has
\begin{equation}
m_{\sigma;\mu,\nu}^{\varepsilon;\varepsilon_1,\varepsilon_2}
=\Phi_{\sigma;\mu,\nu}^{\varepsilon;\varepsilon_1,\varepsilon_2}
\widetilde m_{\sigma;\mu,\nu}^{\varepsilon;\varepsilon_1,\varepsilon_2},
\label{eq:exact-phase-factorization}
\end{equation}
where
\begin{align}
&\widetilde m_{\sigma;\mu,\nu}^{\varepsilon;\varepsilon_1,\varepsilon_2}
(\xi,\eta)=-\frac{i\varepsilon\varepsilon_1}{2\langle\xi\rangle}
\Bigl(\langle\eta\rangle
+\varepsilon_1\varepsilon_2\langle\xi-\eta\rangle
+h\varepsilon_1
+\varepsilon_1\varepsilon\langle\xi\rangle\Bigr).
\label{eq:reduced-normal-form-symbol}
\end{align}
Consequently, $m/\Phi$ admits a smooth extension across the entire
time-resonant set.
\end{proposition}
\begin{proof}
Set
\[
A=\langle\eta\rangle,
\qquad B=\langle\xi-\eta\rangle,
\qquad C=\langle\xi\rangle,
\]
and write $r_1=\varepsilon_1\mu$ and $r_2=\varepsilon_2\nu$.
The discrete coefficient in \eqref{eq:full-bilinear-symbol} is
\[
c=\sigma(r_1+r_2)+\varepsilon(1+r_1r_2).
\]
If $r_1=-r_2$, then $c=0$.
If $r_1=r_2=h$, then $c=2(\sigma h+\varepsilon)$.
Hence $c\neq0$ precisely when
$r_1=r_2=h=\varepsilon\sigma$, which proves
\eqref{eq:active-mode-condition}.
In this case, $c=4\varepsilon$, and the phase reduces to
\begin{equation}
\Phi=h-\varepsilon C+\varepsilon_1A+\varepsilon_2B.
\label{eq:active-phase}
\end{equation}
Moreover,
\[
q_{\mu\nu}(\eta,\xi-\eta)
=\varepsilon_1\varepsilon_2
(A+h\varepsilon_1)(B+h\varepsilon_2)
-\eta\cdot(\xi-\eta).
\]
Using
\[
C^2=A^2+B^2-1+2\eta\cdot(\xi-\eta),
\]
we obtain
\begin{equation}
2q_{\mu\nu}
=\left(A+\varepsilon_1\varepsilon_2B+h\varepsilon_1\right)^2-C^2.
\label{eq:q-difference-squares}
\end{equation}
On the other hand, \eqref{eq:active-phase} gives
\[
A+\varepsilon_1\varepsilon_2B+h\varepsilon_1
-\varepsilon_1\varepsilon C=\varepsilon_1\Phi.
\]
Factoring the difference of squares therefore yields
\begin{equation}
q_{\mu\nu}
=\frac{\varepsilon_1}{2}\Phi
\Bigl(A+\varepsilon_1\varepsilon_2B+h\varepsilon_1
+\varepsilon_1\varepsilon C\Bigr).
\label{eq:q-phase-factorization}
\end{equation}
Substitution into \eqref{eq:full-bilinear-symbol}, with
$c=4\varepsilon$, proves \eqref{eq:exact-phase-factorization}
and \eqref{eq:reduced-normal-form-symbol}.
Since $C=\langle\xi\rangle\geq1$, the extended quotient is smooth
on the whole frequency space.
\end{proof}
\begin{remark}
\label{rem:stronger-than-effective-resonance}
The factorization gives vanishing on the whole time-resonant set,
which is stronger than vanishing only on the space--time resonant set.
For an active interaction, \eqref{eq:active-phase} and
\eqref{eq:reduced-normal-form-symbol} also give
\[
\widetilde m(\xi,\eta)
=-i\left(1+\frac{\varepsilon\Phi(\xi,\eta)}
{2\langle\xi\rangle}\right).
\]
Thus division by $\Phi$ introduces no small-phase singularity,
and the quadratic time resonances can be treated by a single
normal-form transformation without integration by parts in $\eta$.
The mapping properties of the resulting multiplier are established
by the dyadic estimates below.
\end{remark}

\subsection{Quadratic Normal-Form Reduction}
Write $\alpha=(\sigma,\varepsilon)$,
$\beta=(\mu,\varepsilon_1)$, and $\gamma=(\nu,\varepsilon_2)$.
The profile equation is
\begin{equation}
\begin{aligned}
\partial_t\widehat f_\alpha(t,\xi)
={}&\frac{1}{(2\pi)^3}
\sum_{\beta,\gamma}\int_{\mathbb R^3}
e^{it\Phi_{\alpha;\beta,\gamma}(\xi,\eta)}
m_{\alpha;\beta,\gamma}(\xi,\eta)
\\
&\qquad\times
\widehat f_\beta(t,\eta)
\widehat f_\gamma(t,\xi-\eta)\,d\eta.
\end{aligned}
\label{eq:profile-equation-normal-form}
\end{equation}
For an inactive interaction, set $b_{\alpha;\beta,\gamma}=0$.
For an active interaction, define
\begin{equation}
b_{\alpha;\beta,\gamma}(\xi,\eta)
:=\frac{\widetilde m_{\alpha;\beta,\gamma}(\xi,\eta)}{i}.
\label{eq:global-normal-form-symbol}
\end{equation}
This agrees with $m/(i\Phi)$ for $\Phi\neq0$ and defines its smooth
extension across the time-resonant set.
Explicitly, with $h$ as in \eqref{eq:active-mode-condition},
\begin{align}
b_{\alpha;\beta,\gamma}(\xi,\eta)
={}&-\frac{\varepsilon\varepsilon_1}{2\langle\xi\rangle}
\Bigl(\langle\eta\rangle
+\varepsilon_1\varepsilon_2\langle\xi-\eta\rangle
+h\varepsilon_1+\varepsilon_1\varepsilon\langle\xi\rangle\Bigr).
\label{eq:explicit-normal-form-symbol}
\end{align}
Define the quadratic correction by
\begin{align*}
\widehat{\mathfrak B_\alpha[f,f]}(t,\xi)
:={}&\frac{1}{(2\pi)^3}
\sum_{\beta,\gamma}\int_{\mathbb R^3}
e^{it\Phi_{\alpha;\beta,\gamma}(\xi,\eta)}
b_{\alpha;\beta,\gamma}(\xi,\eta)
\\
&\qquad\times
\widehat f_\beta(t,\eta)
\widehat f_\gamma(t,\xi-\eta)\,d\eta,
\end{align*}
and introduce the modified profile
\begin{equation}
F_\alpha:=f_\alpha-\mathfrak B_\alpha[f,f].
\label{eq:modified-profile}
\end{equation}
Since $i\Phi_{\alpha;\beta,\gamma}b_{\alpha;\beta,\gamma}
=m_{\alpha;\beta,\gamma}$, the derivative of the oscillatory factor
cancels the right-hand side of \eqref{eq:profile-equation-normal-form}.
Consequently,
\begin{align}
\partial_t\widehat F_\alpha(t,\xi)
={}&-\frac{1}{(2\pi)^3}
\sum_{\beta,\gamma}\int_{\mathbb R^3}
e^{it\Phi_{\alpha;\beta,\gamma}(\xi,\eta)}
b_{\alpha;\beta,\gamma}(\xi,\eta)
\notag\\
&\quad\times\Bigl(
\partial_t\widehat f_\beta(t,\eta)
\widehat f_\gamma(t,\xi-\eta)
\notag\\
&\hspace{1.8cm}
+\widehat f_\beta(t,\eta)
\partial_t\widehat f_\gamma(t,\xi-\eta)
\Bigr)\,d\eta.
\label{eq:normal-form-generated-profile-equation}
\end{align}
Substituting the quadratic profile equation for each time derivative
makes the right-hand side exactly cubic in the original profiles $f$.
Let $\mathfrak B_{\alpha;\beta,\gamma}$ denote the time-independent
bilinear operator with symbol $b_{\alpha;\beta,\gamma}$ and the
same convolution normalization.
Since $U_\alpha(t)=e^{it\omega_\alpha(D)}f_\alpha(t)$, we have
\[
\mathfrak B_\alpha[f,f](t)
=e^{-it\omega_\alpha(D)}
\sum_{\beta,\gamma}
\mathfrak B_{\alpha;\beta,\gamma}
\bigl(U_\beta(t),U_\gamma(t)\bigr).
\]
For a fixed active interaction, write
$\mathfrak B_b=\mathfrak B_{\alpha;\beta,\gamma}$.
The explicit symbol gives the operator identity
\begin{align}
\mathfrak B_b(u,v)
={}&-\frac{\varepsilon\varepsilon_1}{2}
\langle D\rangle^{-1}
\Bigl((\langle D\rangle u)v
+\varepsilon_1\varepsilon_2u(\langle D\rangle v)
+h\varepsilon_1uv\Bigr)-\frac12uv.
\label{eq:normal-form-operator-decomposition}
\end{align}
The low-frequency gains in the derivative estimates depend on
branch-dependent cancellations among these terms.
The dyadic estimates below retain these cancellations by treating
the combined symbol.

\subsection{Dyadic Multiplier Estimates for the Normal Form}
\begin{lemma}[Derivative-adapted dyadic multiplier bounds]
\label{lem:dyadic-multiplier-bounds}
Let \(\varphi_k\) be the symbol of \(P_k\), set
\(k^+=\max\{k,0\}\), and define
\begin{equation}
d_+(k):=2^{k^+},
\qquad d_-(k):=2^k.
\label{eq:dyadic-modal-weight}
\end{equation}
For every \(1<p<\infty\), and also for frequency-localized
\(L^\infty\) norms,
\begin{equation}
\|\mathcal D_{(\rho,e)}P_kh\|_{L^p}
\simeq d_\rho(k)\|P_kh\|_{L^p}.
\label{eq:D-dyadic-equivalence}
\end{equation}
For an active interaction \((\beta,\gamma)\to\alpha\), where
\(\alpha=(\sigma,\varepsilon)\),
\(\beta=(\mu,\varepsilon_1)\), and
\(\gamma=(\nu,\varepsilon_2)\), let
\begin{align*}
b_{k,k_1,k_2}(\xi,\eta)
&:=\varphi_k(\xi)\varphi_{k_1}(\eta)
\varphi_{k_2}(\xi-\eta)b_{\alpha;\beta,\gamma}(\xi,\eta),
\\
m_{k,k_1,k_2}(\xi,\eta)
&:=\varphi_k(\xi)\varphi_{k_1}(\eta)
\varphi_{k_2}(\xi-\eta)m_{\alpha;\beta,\gamma}(\xi,\eta).
\end{align*}
If \(\mathcal K\) denotes the bilinear kernel, then
\begin{align}
\|\mathcal K[b_{k,k_1,k_2}]\|_{L^1(\mathbb R^6)}
&\lesssim1+d_\mu(k_1)+d_\nu(k_2),
\label{eq:b-dyadic-kernel-bound}
\\
d_\sigma(k)\|\mathcal K[b_{k,k_1,k_2}]\|_{L^1(\mathbb R^6)}
&\lesssim d_\mu(k_1)+d_\nu(k_2),
\label{eq:Db-dyadic-kernel-bound}
\\
\bigl(1+d_\sigma(k)\bigr)
\|\mathcal K[m_{k,k_1,k_2}]\|_{L^1(\mathbb R^6)}
&\lesssim d_\mu(k_1)d_\nu(k_2).
\label{eq:m-dyadic-kernel-bound}
\end{align}
More precisely, for spatial multi-indices \(a,c\) with
\(|a|+|c|\leq2\), set
\begin{equation}
\Theta_{a,c}(k,k_1,k_2)
:=\bigl(2^{-k}+2^{-k_2}\bigr)^{|a|}
\bigl(2^{-k_1}+2^{-k_2}\bigr)^{|c|}.
\label{eq:dyadic-symbol-derivative-factor}
\end{equation}
Then
\begin{align}
\|\mathcal K[\partial_\xi^a\partial_\eta^c
b_{k,k_1,k_2}]\|_{L^1}
&\lesssim
\Theta_{a,c}\bigl(1+d_\mu(k_1)+d_\nu(k_2)\bigr),
\notag\\
d_\sigma(k)\|\mathcal K[\partial_\xi^a\partial_\eta^c
b_{k,k_1,k_2}]\|_{L^1}
&\lesssim
\Theta_{a,c}\bigl(d_\mu(k_1)+d_\nu(k_2)\bigr),
\notag\\
\bigl(1+d_\sigma(k)\bigr)
\|\mathcal K[\partial_\xi^a\partial_\eta^c
m_{k,k_1,k_2}]\|_{L^1}
&\lesssim
\Theta_{a,c}d_\mu(k_1)d_\nu(k_2).
\label{eq:dyadic-differentiated-kernel-bounds}
\end{align}
Here \(\Theta_{a,c}=\Theta_{a,c}(k,k_1,k_2)\), and all kernel
norms are taken over \(\mathbb R^6\). The constants are uniform
in the dyadic indices and the discrete signs.
\end{lemma}

\begin{proof}
Write \(\Lambda(\zeta)=\langle\zeta\rangle\), so that
\(\lambda_\rho(\zeta)=\Lambda(\zeta)+\rho\).
We first record the one-variable kernel estimates used below.
For a one-variable symbol \(q\), write
\(\check q=q^\vee=\mathcal F^{-1}q\).
Rescaling a single annulus by \(\zeta=2^j\theta\) gives,
for every multi-index \(a\) with \(|a|\leq2\),
\begin{align}
\|[\partial^a\varphi_j]^\vee\|_{L^1}
&\lesssim 2^{-j|a|},
\notag\\
\|[\partial^a(\varphi_j\Lambda)]^\vee\|_{L^1}
&\lesssim 2^{-j|a|}2^{j^+},
\notag\\
\|[\partial^a(\varphi_j\Lambda^{-1})]^\vee\|_{L^1}
&\lesssim 2^{-j|a|}2^{-j^+},
\notag\\
\|[\partial^a(\varphi_j(\Lambda-1))]^\vee\|_{L^1}
&\lesssim 2^{-j|a|}2^{2j-j^+},
\notag\\
\|[\partial^a(\varphi_j\lambda_\rho)]^\vee\|_{L^1}
&\lesssim 2^{-j|a|}d_\rho(j),
\notag\\
\|[\partial^a(\varphi_j\zeta_\ell)]^\vee\|_{L^1}
&\lesssim 2^{-j|a|}2^j,
\qquad \ell=1,2,3.
\label{eq:one-variable-dyadic-kernels}
\end{align}
Indeed, after removing the displayed size factors, the rescaled
symbols have uniformly bounded derivatives on a fixed annulus.
For a smooth symbol \(q\) supported there, integration by parts
with \((1-\Delta_\theta)^2\) gives
\[
\|\check q\|_{L^1(\mathbb R^3)}
\lesssim\sum_{|b|\leq4}\|\partial_\theta^bq\|_{L^\infty}.
\]
The additional low-frequency gain in the fourth estimate follows
from \(\Lambda-1=|\zeta|^2/(\Lambda+1)\).
These estimates also hold for slightly enlarged annular cutoffs.

Recall that
\(\mathcal D_{(\rho,e)}h=(\nabla h,\lambda_\rho(D)h)\).
The upper bound in \eqref{eq:D-dyadic-equivalence} follows from
\eqref{eq:one-variable-dyadic-kernels} and Young's inequality.
For the reverse bound, choose a cutoff \(\widetilde\varphi_k\)
equal to one on the support of \(\varphi_k\).
If \(\rho=-1\), or if \(k\geq0\), use
\[
P_kh
=\sum_{\ell=1}^3
\left(\frac{-i\zeta_\ell}{|\zeta|^2}
\widetilde\varphi_k(\zeta)\right)(D)
\partial_\ell P_kh.
\]
Each multiplier on the right has kernel norm \(O(2^{-k})\),
and \(d_\rho(k)=2^k\) in these cases.
If \(\rho=+1\) and \(k<0\), use instead
\[
P_kh
=\left(\frac{\widetilde\varphi_k}{\lambda_+}\right)(D)
\lambda_+(D)P_kh,
\]
whose multiplier has uniformly bounded kernel norm.
This proves the equivalence, including the localized
\(L^\infty\) statement.
To estimate the bilinear kernels, use the independent input
frequencies \((p,q)=(\eta,\xi-\eta)\).
The change from \((\xi,\eta)\) to \((p,q)\) has determinant of
absolute value one and preserves the kernel \(L^1\) norm.
For a separated symbol
\(A(\xi)B(\eta)C(\xi-\eta)\), its kernel in these coordinates is
\[
\int_{\mathbb R^3}
\check A(w)\check B(y-w)\check C(z-w)\,dw.
\]
Consequently, Fubini's theorem gives
\begin{equation}
\|\mathcal K[A(\xi)B(\eta)C(\xi-\eta)]\|_{L^1(\mathbb R^6)}
\leq\|\check A\|_{L^1}
\|\check B\|_{L^1}\|\check C\|_{L^1}.
\label{eq:separated-bilinear-kernel}
\end{equation}
This estimate allows all three cutoffs to retain their own scales.
Set
\[
r=2^k,\quad r_1=2^{k_1},\quad r_2=2^{k_2},
\qquad
L=2^{k^+},\quad L_1=2^{k_1^+},\quad L_2=2^{k_2^+}.
\]
If the product of the cutoffs is nonzero, then
\begin{equation}
r\lesssim r_1+r_2,
\qquad L\lesssim L_1+L_2.
\label{eq:dyadic-output-triangle}
\end{equation}
If that product vanishes identically, the desired bounds are
immediate.
For an active interaction,
\(\varepsilon\varepsilon_1=\sigma\mu\),
\(\varepsilon\varepsilon_2=\sigma\nu\), and
\(\varepsilon h=\sigma\).
Thus \eqref{eq:explicit-normal-form-symbol} becomes
\begin{equation}
b_{\alpha;\beta,\gamma}(\xi,\eta)
=-\frac12
-\frac{\sigma}{2\Lambda(\xi)}
\bigl(\mu\Lambda(\eta)+\nu\Lambda(\xi-\eta)+1\bigr).
\label{eq:branchwise-normal-form-symbol}
\end{equation}
Each term, after localization, has the separated form in
\eqref{eq:separated-bilinear-kernel}. Hence
\begin{equation}
\|\mathcal K[b_{k,k_1,k_2}]\|_{L^1}
\lesssim1+\frac{L_1+L_2+1}{L}
\lesssim\frac{L_1+L_2}{L}.
\label{eq:b-kernel-preliminary}
\end{equation}
Since \(L\geq1\) and
\(L_1+L_2\lesssim1+d_\mu(k_1)+d_\nu(k_2)\), this proves
\eqref{eq:b-dyadic-kernel-bound}.

For \eqref{eq:Db-dyadic-kernel-bound}, suppose first that at
least one input is a fast mode, or that
\(\max\{k_1,k_2\}\geq0\). Then
\[
L_1+L_2\lesssim d_\mu(k_1)+d_\nu(k_2).
\]
Since \(d_\sigma(k)\leq L\),
\eqref{eq:b-kernel-preliminary} gives the desired estimate.

It remains to consider two slow inputs at low frequency:
\(\mu=\nu=-1\) and \(k_1,k_2<0\).
If the output is slow, then \(d_-(k)=r\), and
\eqref{eq:b-kernel-preliminary} gives
\[
d_-(k)\|\mathcal K[b_{k,k_1,k_2}]\|_{L^1}
\lesssim r\lesssim r_1+r_2
=d_-(k_1)+d_-(k_2).
\]
If the output is fast, the constant terms in
\eqref{eq:branchwise-normal-form-symbol} cancel exactly:
\begin{equation}
b_{\alpha;\beta,\gamma}(\xi,\eta)
=\frac{
(\Lambda(\eta)-1)+(\Lambda(\xi-\eta)-1)
-(\Lambda(\xi)-1)}{2\Lambda(\xi)}.
\label{eq:slow-slow-fast-b-cancellation}
\end{equation}
Use the fourth estimate in
\eqref{eq:one-variable-dyadic-kernels}, together with
\eqref{eq:separated-bilinear-kernel}, to obtain
\[
\|\mathcal K[b_{k,k_1,k_2}]\|_{L^1}
\lesssim\frac{r_1^2+r_2^2}{L}+\frac{r^2}{L^2}.
\]
Here the output multiplier
\(\varphi_k(\Lambda-1)/\Lambda\) has kernel norm
\(O(r^2/L^2)\), by the same one-variable rescaling.
Since \(d_+(k)=L\), \(r_j<1\) for \(j=1,2\), and
\(r^2/L\leq r\), we conclude that
\[
d_+(k)\|\mathcal K[b_{k,k_1,k_2}]\|_{L^1}
\lesssim r_1^2+r_2^2+\frac{r^2}{L}
\lesssim r_1+r_2.
\]
This completes the proof of
\eqref{eq:Db-dyadic-kernel-bound}.
For the quadratic symbol, the active relations reduce
\eqref{eq:full-bilinear-symbol} to
\[
m_{\alpha;\beta,\gamma}(\xi,\eta)
=-\frac{i\varepsilon}{\Lambda(\xi)}
\Bigl(
\varepsilon_1\varepsilon_2
\lambda_\mu(\eta)\lambda_\nu(\xi-\eta)
-\eta\cdot(\xi-\eta)
\Bigr).
\]
Applying \eqref{eq:one-variable-dyadic-kernels} and
\eqref{eq:separated-bilinear-kernel} term by term gives
\[
\|\mathcal K[m_{k,k_1,k_2}]\|_{L^1}
\lesssim L^{-1}
\bigl(d_\mu(k_1)d_\nu(k_2)+r_1r_2\bigr)
\lesssim L^{-1}d_\mu(k_1)d_\nu(k_2).
\]
Because \(1+d_\sigma(k)\leq2L\), this proves
\eqref{eq:m-dyadic-kernel-bound}.
Finally, differentiate the fully localized separated symbols.
A derivative \(\partial_\xi\) acts on the output factor or
on the factor depending on \(\xi-\eta\), and therefore costs
at most \(2^{-k}+2^{-k_2}\).
A derivative \(\partial_\eta\) acts on the factor depending
on \(\eta\) or on the factor depending on \(\xi-\eta\),
and costs at most \(2^{-k_1}+2^{-k_2}\).
The Leibniz rule and the differentiated estimates in
\eqref{eq:one-variable-dyadic-kernels} give the factor
\(\Theta_{a,c}\).
In the slow--slow to fast case, differentiate the already
cancelled expression \eqref{eq:slow-slow-fast-b-cancellation};
the localized bounds for \(\Lambda-1\) retain the gains used
above, with the same inverse dyadic factors.
Thus every preceding kernel estimate remains valid after
\(\partial_\xi^a\partial_\eta^c\), with the additional factor
\(\Theta_{a,c}\). This proves
\eqref{eq:dyadic-differentiated-kernel-bounds}.
\end{proof}

For a scalar function \(h\), define
\begin{equation}
\|h\|_{\mathcal B_K}
:=
\|P_{\leq0}h\|_{L^\infty}
+
\sum_{k\geq1}2^{Kk}\|P_kh\|_{L^\infty}.
\label{eq:BK-definition}
\end{equation}
\begin{corollary}[Tame bounds for the quadratic normal form]
\label{cor:normal-form-tame-bounds}
Let \(s\geq0\) and \(r\geq K+3\). Use the normalized
representatives for slow modes fixed above.
For every active interaction \((\beta,\gamma)\to\alpha\),
the following estimates hold whenever their right-hand sides
are finite:
\begin{align}
\|\mathcal D_\alpha
\mathfrak B_{\alpha;\beta,\gamma}(u,v)\|_{H^s}
\lesssim{}&
\|\mathcal D_\beta u\|_{H^s}
\|\mathcal D_\gamma v\|_{H^r}
+\|\mathcal D_\beta u\|_{H^r}
\|\mathcal D_\gamma v\|_{H^s},
\label{eq:normal-form-derivative-tame-Hs}
\\
\|\mathcal D_\alpha
\mathfrak B_{\alpha;\beta,\gamma}(u,v)\|_{\mathcal B_K}
\lesssim{}&
\|\mathcal D_\beta u\|_{\mathcal B_K}
\|\mathcal D_\gamma v\|_{H^r}
+\|\mathcal D_\beta u\|_{H^r}
\|\mathcal D_\gamma v\|_{\mathcal B_K}.
\label{eq:normal-form-derivative-tame-BK}
\end{align}
When the lower-order factors are measured in \(\mathcal B_K\),
we have
\begin{align}
\|\mathcal D_\alpha
\mathfrak B_{\alpha;\beta,\gamma}(u,v)\|_{H^s}
\lesssim{}&
\|\mathcal D_\beta u\|_{H^s}
\Bigl(\|v\|_{\mathcal B_K}
+\|\mathcal D_\gamma v\|_{\mathcal B_K}\Bigr)
\notag\\
&+\Bigl(\|u\|_{\mathcal B_K}
+\|\mathcal D_\beta u\|_{\mathcal B_K}\Bigr)
\|\mathcal D_\gamma v\|_{H^s},
\label{eq:normal-form-pointwise-tame-Hs}
\\
\|\mathfrak B_{\alpha;\beta,\gamma}(u,v)\|_{\mathcal B_K}
+\|\mathcal D_\alpha
\mathfrak B_{\alpha;\beta,\gamma}(u,v)\|_{\mathcal B_K}
\lesssim{}&
\Bigl(\|u\|_{\mathcal B_K}
+\|\mathcal D_\beta u\|_{\mathcal B_K}\Bigr)
\Bigl(\|v\|_{\mathcal B_K}
+\|\mathcal D_\gamma v\|_{\mathcal B_K}\Bigr).
\label{eq:normal-form-pointwise-tame-BK}
\end{align}
The constants may depend on \(s,r,K\) and the fixed cutoffs,
but are uniform in the discrete indices.
Every Sobolev norm on the right-hand sides is applied to a modal
derivative. No independent \(L^2\) norm of an undifferentiated
slow input is required.
\end{corollary}

\begin{proof}
We first explain how the derivative Sobolev norm controls the
pointwise norms of a normalized input. In three dimensions,
Bernstein's inequality and Cauchy--Schwarz give
\begin{align}
\|P_{\leq0}h\|_{L^\infty}
&\lesssim\sum_{j\leq0}2^{3j/2}\|P_jh\|_{L^2}
\notag\\
&\lesssim\sum_{j\leq0}2^{j/2}\|P_j\nabla h\|_{L^2}
\lesssim\|\nabla h\|_{L^2}.
\label{eq:tame-low-frequency-pointwise}
\end{align}
The normalization removes the constant mode, so that the
low-frequency reconstruction used here is valid.
At high frequency,
\[
\sum_{j\geq1}2^{Kj}\|P_jh\|_{L^\infty}
\lesssim
\sum_{j\geq1}2^{(K+1/2-r)j}
\bigl(2^{rj}\|P_j\nabla h\|_{L^2}\bigr)
\lesssim\|\nabla h\|_{H^r}.
\]
Also, \(H^r(\mathbb R^3)\hookrightarrow\mathcal B_K\)
because \(r>K+3/2\). Hence, for every modal index \(\theta\),
\begin{equation}
\|h\|_{\mathcal B_K}
+\|\mathcal D_\theta h\|_{\mathcal B_K}
\lesssim\|\mathcal D_\theta h\|_{H^r}.
\label{eq:modal-derivative-Sobolev-to-BK}
\end{equation}
We next prove the estimates with pointwise lower-order factors.
Write \(\Lambda=\langle D\rangle\) and
\(\lambda_-(D)=\Lambda-1\).
The explicit symbol in
\eqref{eq:branchwise-normal-form-symbol} gives
\begin{align}
\mathfrak B_{\alpha;\beta,\gamma}(u,v)
=-\frac12\Lambda^{-1}\Bigl[
&c_{\sigma\mu\nu}uv+\lambda_-(D)(uv)
\notag\\
&+\sigma\mu\bigl(\lambda_-(D)u\bigr)v
+\sigma\nu u\bigl(\lambda_-(D)v\bigr)
\Bigr],
\label{eq:tame-normal-form-decomposition}
\end{align}
where
\[
c_{\sigma\mu\nu}:=1+\sigma(1+\mu+\nu).
\]
The localized multiplier estimates in
Lemma~\ref{lem:dyadic-multiplier-bounds} imply that
\(\Lambda^{-1}\), \(\lambda_-(D)\Lambda^{-1}\), and
\(\mathcal D_\alpha\Lambda^{-1}\) are bounded on both
\(H^s\) and \(\mathcal B_K\).
Moreover,
\begin{equation}
\|\lambda_-(D)h\|_X\lesssim\|\nabla h\|_X,
\qquad X=H^s\ \hbox{or}\ \mathcal B_K,
\label{eq:lambda-minus-by-gradient-tame}
\end{equation}
since
\[
\lambda_-(D)h
=\sum_{\ell=1}^3T_\ell(D)\partial_\ell h,
\qquad
T_\ell(\xi)=\frac{-i\xi_\ell}{\langle\xi\rangle+1},
\]
and the \(T_\ell\) are bounded order-zero multipliers on these
spaces.
The constant term in
\eqref{eq:tame-normal-form-decomposition} is handled as follows.
If \(\sigma=-1\), then
\[
\|\mathcal D_\alpha\Lambda^{-1}(uv)\|_X
\lesssim\|\nabla(uv)\|_X.
\]
If \(\sigma=+1\) and \(\mu=\nu=-1\), then
\(c_{\sigma\mu\nu}=0\).
In the remaining cases the output is fast and at least one input
is fast. That input itself is controlled by its modal derivative,
because \(\lambda_+(D)^{-1}\) is bounded on \(X\).
Thus the constant term also admits the derivative estimates below.
The algebra property of \(\mathcal B_K\), applied to
\eqref{eq:tame-normal-form-decomposition}, now gives the useful
intermediate estimate
\begin{align}
\|\mathcal D_\alpha
\mathfrak B_{\alpha;\beta,\gamma}(u,v)\|_{\mathcal B_K}
\lesssim{}&
\|\mathcal D_\beta u\|_{\mathcal B_K}\|v\|_{\mathcal B_K}
\notag\\
&+\|u\|_{\mathcal B_K}
\|\mathcal D_\gamma v\|_{\mathcal B_K}.
\label{eq:normal-form-pointwise-derivative-BK}
\end{align}
Here \(\lambda_-(D)(uv)\) is estimated through
\eqref{eq:lambda-minus-by-gradient-tame} and the ordinary
product rule for \(\nabla(uv)\).
Without the output derivative, the same decomposition gives
\begin{align*}
\|\mathfrak B_{\alpha;\beta,\gamma}(u,v)\|_{\mathcal B_K}
\lesssim{}&\|u\|_{\mathcal B_K}\|v\|_{\mathcal B_K}
+\|\mathcal D_\beta u\|_{\mathcal B_K}\|v\|_{\mathcal B_K}
\\
&+\|u\|_{\mathcal B_K}
\|\mathcal D_\gamma v\|_{\mathcal B_K}.
\end{align*}
These two bounds imply
\eqref{eq:normal-form-pointwise-tame-BK}.
For the Sobolev estimate, we use the product inequality
\begin{equation}
\|fg\|_{H^s}
\lesssim_s
\|f\|_{H^s}\|g\|_{\mathcal B_K}
+\|f\|_{\mathcal B_K}\|\nabla g\|_{H^s}.
\label{eq:mixed-product-no-undifferentiated-L2}
\end{equation}
To see this, split \(g=P_{\leq0}g+(1-P_{\leq0})g\).
Multiplication by \(P_{\leq0}g\) is bounded on \(H^s\),
with norm controlled by \(\|g\|_{L^\infty}\), since all its
spatial derivatives are bounded by the same norm.
For the high-frequency part, use the usual tame Sobolev
product estimate and
\[
\|(1-P_{\leq0})g\|_{H^s}
\lesssim\|\nabla g\|_{H^s}.
\]
For \(s=0\), the required product bound follows directly from
H\"older's inequality.
Apply \eqref{eq:mixed-product-no-undifferentiated-L2} to the
products in \eqref{eq:tame-normal-form-decomposition}, after
applying \(\mathcal D_\alpha\).
For example,
\begin{align*}
\|\bigl(\lambda_-(D)u\bigr)v\|_{H^s}
\lesssim{}&
\|\mathcal D_\beta u\|_{H^s}\|v\|_{\mathcal B_K}
\\
&+\|\mathcal D_\beta u\|_{\mathcal B_K}
\|\mathcal D_\gamma v\|_{H^s}.
\end{align*}
The terms from \(\nabla(uv)\) satisfy the same type of
estimate, and the constant term is treated by the cases
described above. Interchanging the inputs where necessary
proves \eqref{eq:normal-form-pointwise-tame-Hs}.
Finally, substitute \eqref{eq:modal-derivative-Sobolev-to-BK}
into \eqref{eq:normal-form-pointwise-tame-Hs} to obtain
\eqref{eq:normal-form-derivative-tame-Hs}.
Substituting the same embedding into
\eqref{eq:normal-form-pointwise-derivative-BK} gives
\eqref{eq:normal-form-derivative-tame-BK}.
\end{proof}

\section{Small-data global existence}
\label{sec:small-data-global}

We now close the global argument.  The proof follows the normal-form
strategy for three-dimensional quadratic dispersive systems: the
quadratic part is removed first, the cubic Duhamel term generated by
that change of variable is estimated by separating early and late
interaction times, and the
integrable pointwise decay is then inserted into the high-order energy
estimate.  This is the same general mechanism used in
\cite{OzawaTsutayaTsutsumi1995,Shatah1985}, but here all four shifted
branches and the low-frequency Schr\"odinger regime are retained. To avoid confusing the sign index
\(\varepsilon\in\{\pm1\}\) with the small parameter in
\eqref{eq:main-data}, we denote the actual size of the initial data by
\(\delta_0\).

\subsection{Norms and statement of the theorem}
We label the four linear modes by
\begin{equation}
\mathcal A:=\{(\sigma,\varepsilon):
\sigma,\varepsilon\in\{\pm1\}\},
\label{eq:modal-index-set}
\end{equation}
where \(\sigma\) distinguishes the fast and slow branches,
and \(\varepsilon\) specifies the sign of the corresponding
temporal frequency.
Let \(P_k\), \(k\in\mathbb Z\), be homogeneous
Littlewood--Paley projections localizing to frequencies
\(|\xi|\sim2^k\). We denote the low-frequency projection by
\[
P_{\leq0}:=\sum_{k\leq0}P_k,
\]
which localizes to \(|\xi|\lesssim1\).
Fix
\begin{equation}
K=6,
\qquad
M=K+4=10,
\qquad
N\geq30.
\label{eq:KNM-choice}
\end{equation}
The first term controls the low-frequency part, while the
high-frequency weights account for \(K\) spatial derivatives.
For \(K\geq1\), the space \(\mathcal B_K\) is an algebra and
embeds continuously into \(W^{K-1,\infty}\).
For vector-valued functions, we use the sum of the component norms.
With the additive constants of \(\psi\) and \(S\) fixed,
we define the high-order energy norm
\begin{equation}
\begin{aligned}
E_N(t)
:={}&
\|\psi_t(t)\|_{H^N}
+\|\nabla\psi(t)\|_{H^N}+
\|S_t(t)\|_{H^N}
+\|\nabla S(t)\|_{H^N}.
\end{aligned}
\label{eq:EN-modal}
\end{equation}
To measure the pointwise size and spatial regularity of the
first derivatives, we set
\begin{equation}
\begin{aligned}
Z_K(t)
:={}&
\|\psi_t(t)\|_{\mathcal B_K}
+\|\nabla\psi(t)\|_{\mathcal B_K}+
\|S_t(t)\|_{\mathcal B_K}
+\|\nabla S(t)\|_{\mathcal B_K}.
\end{aligned}
\label{eq:ZK-modal}
\end{equation}
We also record the pointwise norms of the unknowns themselves:
\begin{equation}
Y_K(t)
:=
\|\psi(t)\|_{\mathcal B_K}
+
\|S(t)\|_{\mathcal B_K}.
\label{eq:YK-modal}
\end{equation}
Thus, \(E_N\) controls the high-order derivative energy,
\(Z_K\) controls the pointwise derivative norms, and
\(Y_K\) controls the pointwise norms of \(\psi\) and \(S\).
No undifferentiated \(L^2\) norm of either unknown is included
in \(E_N\).
The auxiliary norm \(Y_K\) is included to track undifferentiated
factors arising in the nonlinear estimates after the normal-form
transformation. By the definition of \(\mathcal B_K\), we also have
\begin{equation}
\|\psi_t(t)\|_{L^\infty}
+\|\nabla\psi(t)\|_{L^\infty}
+\|S_t(t)\|_{L^\infty}
+\|\nabla S(t)\|_{L^\infty}
\lesssim Z_K(t).
\label{eq:first-derivatives-by-ZK}
\end{equation}
For a multi-index
\(\rho=(\rho_1,\rho_2,\rho_3)\in\mathbb N_0^3\), write
\[
|\rho|=\rho_1+\rho_2+\rho_3,
\qquad
x^\rho=x_1^{\rho_1}x_2^{\rho_2}x_3^{\rho_3}.
\]
We fix the additive constants of \(\psi_0\) and \(S_0\)
by choosing their representatives that vanish at spatial infinity.
Define the size of the initial profiles by
\begin{align}
\mathcal I_N^\partial
:={}&
\|\nabla\psi_0\|_{H^N}
+\|\nabla S_0\|_{H^N}
+\|\psi_1\|_{H^N}
+\|S_1\|_{H^N}
\notag\\
&+
\sum_{|\rho|\leq3}
\Big(
\|x^\rho\nabla\psi_0\|_{H^M}
+\|x^\rho\nabla S_0\|_{H^M}
\notag\\
&\hspace{2.8cm}
+\|x^\rho\psi_1\|_{H^M}
+\|x^\rho S_1\|_{H^M}
\Big).
\label{eq:physical-initial-size}
\end{align}
For the scaled initial data in \eqref{eq:main-data}, define
\begin{equation}
\delta_0:=\delta\,\mathcal I_N^\partial.
\label{eq:delta0-definition}
\end{equation}
This quantity measures the size of the actual initial data
in the high-order derivative and weighted derivative norms.
No independent \(L^2\) assumption on \(\psi_0\) or \(S_0\)
is required.

\begin{theorem}[Global existence, decay, and scattering]
\label{thm:small-data-global}
Fix \(K=6\), \(M=10\), and an integer \(N\geq30\).
Assume that the real-valued initial profiles are normalized
as above and satisfy \(\mathcal I_N^\partial<\infty\).
There exists
\(\delta_*=\delta_*(\mathcal I_N^\partial)>0\)
such that, for \(0<\delta\leq\delta_*\), the Cauchy problem
\eqref{eq:main-system}--\eqref{eq:main-data}
admits a unique global classical solution satisfying
\begin{equation}
\sup_{t\geq0}E_N(t)
+
\sup_{t\geq0}\langle t\rangle^{3/2}
\bigl(Z_K(t)+Y_K(t)\bigr)
\leq C\delta_0.
\label{eq:main-global-bound}
\end{equation}
Here \(E_N\), \(Z_K\), and \(Y_K\) are the norms of
\(\psi\) and \(S\) defined above. In particular,
\begin{align}
&\|\psi_t(t)\|_{W^{K-1,\infty}}
+\|\nabla\psi(t)\|_{W^{K-1,\infty}}
\notag\\
&\qquad
+\|S_t(t)\|_{W^{K-1,\infty}}
+\|\nabla S(t)\|_{W^{K-1,\infty}}
\lesssim
\delta_0\langle t\rangle^{-3/2}.
\label{eq:main-derivative-decay}
\end{align}
Moreover, there exists a real-valued solution
\((\widetilde{\psi},\widetilde{S})\)
of the linearized system \eqref{eq:linearized-system} such that
\begin{equation}
\begin{aligned}
\lim_{t\to\infty}\Big(
&\|\partial_t(\psi-\widetilde{\psi})(t)\|_{H^N}
+\|\nabla(\psi-\widetilde{\psi})(t)\|_{H^N}\\
&+\|\partial_t(S-\widetilde{S})(t)\|_{H^N}
+\|\nabla(S-\widetilde{S})(t)\|_{H^N}
\Big)=0.
\end{aligned}
\label{eq:physical-scattering}
\end{equation}
\end{theorem}

\begin{remark}[Weighted control of the undifferentiated data]
The third spatial weight allows us to recover two spatial
weights of \(\psi_0\) and \(S_0\) from their derivative data.
Indeed, for \(h\in C_c^\infty(\mathbb R^3)\),
integration by parts gives
\[
(3+2j)\int_{\mathbb R^3}|x|^{2j}|h|^2\,dx
=
-2\operatorname{Re}
\int_{\mathbb R^3}
|x|^{2j}\overline h\,x\cdot\nabla h\,dx.
\]
Hence
\begin{equation}
\||x|^j h\|_{L^2}
\leq
\frac{2}{3+2j}
\||x|^{j+1}\nabla h\|_{L^2},
\qquad j=0,1,2.
\label{eq:weighted-hardy-position-from-derivative}
\end{equation}
Applying this \(L^2\) estimate to spatial derivatives of \(h\)
and using the Leibniz rule yields
\begin{equation}
\sum_{|\rho|\leq2}\|x^\rho h\|_{H^M}
\lesssim_M
\sum_{|\rho|\leq3}\|x^\rho\nabla h\|_{H^M}.
\label{eq:weighted-data-from-derivatives}
\end{equation}
By approximation, the estimate also applies to the normalized
initial profiles with finite right-hand side.
Consequently, their undifferentiated weighted \(H^M\) norms
follow from \(\mathcal I_N^\partial<\infty\), rather than
being imposed as separate assumptions.
\end{remark}

\subsection{Local theory and the high-order energy estimate}

\begin{proposition}[Local existence and continuation criterion]
\label{prop:local-continuation}
Let \(N\geq3\), and assume that
\[
\nabla\psi(0),\quad \nabla S(0),\quad
\partial_t\psi(0),\quad \partial_tS(0)\in H^N(\mathbb R^3).
\]
After fixing the additive constants of the initial potentials,
there is a unique solution of \eqref{eq:main-system}
on a maximal interval \([0,T_*)\) such that
\[
\nabla\psi,\ \nabla S,\ \partial_t\psi,\ \partial_tS
\in C([0,T_*);H^N(\mathbb R^3)).
\]
If \(T_*<\infty\), then
\begin{equation}
\limsup_{t\to T_*}E_N(t)=\infty.
\label{eq:continuation-criterion}
\end{equation}
\end{proposition}

\begin{proof}
The linear evolution associated with \eqref{eq:linearized-system}
is uniformly bounded in the derivative norm \(E_N\): its two
coupling terms cancel in the quadratic energy identity.
The corresponding estimate for the forced linear system bounds
the derivative norm by the initial norm and the time integral of
the \(H^N\) norms of the two forcing terms.
The nonlinearities in \eqref{eq:main-system} are quadratic
polynomials of \(\psi_t,\nabla\psi,S_t,\nabla S\).
Since \(H^N(\mathbb R^3)\) is an algebra, these nonlinearities
are locally Lipschitz from the derivative energy space into
\(H^N\times H^N\). Picard iteration for the Duhamel formula
therefore yields the stated local solution, with an existence
time depending only on the size of \(E_N(0)\).
If \(E_N(t)\) remained bounded as \(t\to T_*<\infty\), the same local
construction could be restarted at a time \(t_0<T_*\), with a lifespan
bounded below uniformly in \(t_0\).
Choosing \(t_0\) sufficiently close to \(T_*\) would extend the
solution beyond \(T_*\), contradicting maximality.
This proves \eqref{eq:continuation-criterion}.
\end{proof}

We estimate the first derivatives of \(\psi\) and \(S\) directly
from \eqref{eq:main-system}. To differentiate the energy in time,
we use the quadratic quantity
\begin{equation}
\mathcal E_N(t)
:=\frac12\sum_{|a|\leq N}
\Big(
\|\partial_x^a\psi_t(t)\|_{L^2}^2
+\|\nabla\partial_x^a\psi(t)\|_{L^2}^2
+\|\partial_x^aS_t(t)\|_{L^2}^2
+\|\nabla\partial_x^aS(t)\|_{L^2}^2
\Big).
\label{eq:homogeneous-high-order-energy}
\end{equation}
By the definition of \(E_N\) in \eqref{eq:EN-modal},
\begin{equation}
\mathcal E_N(t)^{1/2}\simeq E_N(t).
\label{eq:homogeneous-energy-modal-equivalence}
\end{equation}
Thus, this quadratic energy measures precisely the same first
derivatives as \(E_N\).

\begin{lemma}[High-order energy estimate]
\label{lem:high-order-energy}
Let \(N\geq3\) be an integer, and let \((\psi,S)\) be the solution
of \eqref{eq:main-system} given by
Proposition~\ref{prop:local-continuation} on \([0,T_*)\).
Then, for almost every \(t\in[0,T_*)\),
\begin{equation}
\left|\frac{d}{dt}\mathcal E_N(t)\right|
\leq C_N Z_K(t)\mathcal E_N(t).
\label{eq:physical-homogeneous-energy-estimate}
\end{equation}
Consequently, for every \(0\leq t<T_*\),
\begin{equation}
E_N(t)
\leq C_N E_N(0)
\exp\left(C_N\int_0^t Z_K(s)\,ds\right).
\label{eq:EN-gronwall-ZK}
\end{equation}
The constants are independent of \(t\) and the solution.
\end{lemma}

\begin{proof}
Denote the quadratic right-hand sides of \eqref{eq:main-system} by
\begin{align*}
F_\psi
&:=S_t^2-|\nabla S|^2-\psi_t^2+|\nabla\psi|^2,\\
F_S
&:=2\nabla\psi\cdot\nabla S-2\psi_tS_t.
\end{align*}
For each spatial multi-index \(a\) with \(|a|\leq N\),
the differentiated equations are
\begin{align*}
\partial_t^2\partial_x^a\psi
-\Delta\partial_x^a\psi
+2\partial_x^aS_t&=\partial_x^aF_\psi,\\
\partial_t^2\partial_x^aS
-\Delta\partial_x^aS
-2\partial_x^a\psi_t&=\partial_x^aF_S.
\end{align*}
Multiply the first equation by \(\partial_x^a\psi_t\),
the second by \(\partial_x^aS_t\), and integrate over
\(\mathbb R^3\). Integration by parts gives, for example,
\[
-\int_{\mathbb R^3}
\Delta\partial_x^a\psi\,\partial_x^a\psi_t\,dx
=\frac12\frac{d}{dt}
\|\nabla\partial_x^a\psi\|_{L^2}^2.
\]
Upon adding the two identities, the linear coupling terms cancel:
\[
-2\int_{\mathbb R^3}
\partial_x^aS_t\,\partial_x^a\psi_t\,dx
+2\int_{\mathbb R^3}
\partial_x^a\psi_t\,\partial_x^aS_t\,dx=0.
\]
Summing over \(|a|\leq N\), we obtain the exact identity
\begin{equation}
\frac{d}{dt}\mathcal E_N(t)
=\sum_{|a|\leq N}\int_{\mathbb R^3}
\Big(
\partial_x^a\psi_t\,\partial_x^aF_\psi
+\partial_x^aS_t\,\partial_x^aF_S
\Big)\,dx.
\label{eq:high-order-energy-identity}
\end{equation}
To estimate the right-hand side, we use the tame product inequality
\begin{equation}
\|fg\|_{H^N}
\leq C_N\Big(
\|f\|_{L^\infty}\|g\|_{H^N}
+\|g\|_{L^\infty}\|f\|_{H^N}
\Big).
\label{eq:energy-tame-product}
\end{equation}
For completeness, the Leibniz rule reads
\[
\partial_x^a(fg)
=\sum_{b\leq a}\binom{a}{b}
\partial_x^bf\,\partial_x^{a-b}g.
\]
The terms with \(b=0\) or \(b=a\) are bounded directly by
the right-hand side of \eqref{eq:energy-tame-product}.
For an intermediate term, set \(m=|a|\) and \(j=|b|\),
where \(1\leq j\leq m-1\). H\"older's inequality and
Gagliardo--Nirenberg interpolation give
\begin{align*}
\|\partial_x^bf\,\partial_x^{a-b}g\|_{L^2}
&\leq
\|\partial_x^bf\|_{L^{2m/j}}
\|\partial_x^{a-b}g\|_{L^{2m/(m-j)}}\\
&\leq C_N
\|f\|_{L^\infty}^{1-j/m}\|f\|_{H^N}^{j/m}
\|g\|_{L^\infty}^{j/m}\|g\|_{H^N}^{1-j/m}\\
&=C_N
\big(\|f\|_{L^\infty}\|g\|_{H^N}\big)^{1-j/m}
\big(\|g\|_{L^\infty}\|f\|_{H^N}\big)^{j/m}\\
&\leq C_N\Big(
\|f\|_{L^\infty}\|g\|_{H^N}
+\|g\|_{L^\infty}\|f\|_{H^N}
\Big).
\end{align*}
The last step is Young's inequality. Summing the finitely many
terms proves \eqref{eq:energy-tame-product}.
In particular, this estimate uses only the \(L^\infty\) norms of
\(f,g\), even when both factors in a Leibniz term carry derivatives.
Applying \eqref{eq:energy-tame-product} to the four square terms
in \(F_\psi\) yields
\begin{align*}
\|F_\psi\|_{H^N}
\leq C_N\Big(&
\|S_t\|_{L^\infty}\|S_t\|_{H^N}
+\|\nabla S\|_{L^\infty}\|\nabla S\|_{H^N}\\
&+\|\psi_t\|_{L^\infty}\|\psi_t\|_{H^N}
+\|\nabla\psi\|_{L^\infty}\|\nabla\psi\|_{H^N}
\Big).
\end{align*}
For the two mixed products in \(F_S\), the same inequality gives
\begin{align*}
\|F_S\|_{H^N}
\leq C_N\Big(&
\|\nabla\psi\|_{L^\infty}\|\nabla S\|_{H^N}
+\|\nabla S\|_{L^\infty}\|\nabla\psi\|_{H^N}\\
&+\|\psi_t\|_{L^\infty}\|S_t\|_{H^N}
+\|S_t\|_{L^\infty}\|\psi_t\|_{H^N}
\Big).
\end{align*}
Since \(\|h\|_{L^\infty}\lesssim\|h\|_{\mathcal B_K}\),
the definitions of \(E_N\) and \(Z_K\) imply
\begin{equation}
\|F_\psi(t)\|_{H^N}+\|F_S(t)\|_{H^N}
\leq C_N Z_K(t)E_N(t).
\label{eq:quadratic-source-energy-bound}
\end{equation}
Thus, every factor measured in \(H^N\) is a first derivative
of \(\psi\) or \(S\).
By Cauchy--Schwarz in \eqref{eq:high-order-energy-identity},
followed by \eqref{eq:quadratic-source-energy-bound},
\begin{align*}
\left|\frac{d}{dt}\mathcal E_N(t)\right|
&\leq C_N E_N(t)
\big(\|F_\psi(t)\|_{H^N}+\|F_S(t)\|_{H^N}\big)\\
&\leq C_N Z_K(t)E_N(t)^2\\
&\leq C_N Z_K(t)\mathcal E_N(t).
\end{align*}
This proves \eqref{eq:physical-homogeneous-energy-estimate}.
Applying the inequality to
\((\mathcal E_N+\epsilon)^{1/2}\) and letting
\(\epsilon\downarrow0\), we obtain
\[
\mathcal E_N(t)^{1/2}
\leq\mathcal E_N(0)^{1/2}
+C_N\int_0^t Z_K(s)\mathcal E_N(s)^{1/2}\,ds.
\]
In particular, the norm equivalence gives
\begin{equation}
E_N(t)\leq C_N E_N(0)
+C_N\int_0^t Z_K(s)E_N(s)\,ds.
\label{eq:EN-integral-estimate}
\end{equation}
Gronwall's inequality proves \eqref{eq:EN-gronwall-ZK}.
The calculation is first performed for smooth solutions;
regularization and passage to the limit give the same estimate
for the solutions of Proposition~\ref{prop:local-continuation}.
\end{proof}
The differential estimate is stated for \(\mathcal E_N\), where
the linear coupling cancels exactly. The equivalence with the sum
norm \(E_N\) then gives the integral and exponential bounds above.
No undifferentiated \(L^2\) norm of \(\psi\) or \(S\) occurs.
In particular, a bound for
\(\int_0^{T_*}Z_K(s)\,ds\) yields uniform control of \(E_N\)
on the lifespan of the solution. By
Proposition~\ref{prop:local-continuation}, such a bound rules out
a finite maximal lifespan. The decay estimates established below
will provide this time-integrability.

\subsection{Linear decay and normal-form-generated cubic estimates}
The dyadic dispersive estimates are naturally stated in terms of the
$L^1$-norm of the initial profile. The weighted $L^2$ assumptions
provide this $L^1$-control and therefore yield the
$\langle t\rangle^{-3/2}$ decay required in the bootstrap.
\begin{lemma}[Weighted linear decay]
\label{lem:weighted-linear-decay}
For \(K,M\) as in \eqref{eq:KNM-choice}, every
\(\alpha\in\mathcal A\), and every \(t\geq0\),
\begin{equation}
\|e^{it\omega_\alpha(D)}h\|_{\mathcal B_K}
\lesssim
\langle t\rangle^{-3/2}\|h\|_{W^{M,1}}.
\label{eq:WM1-linear-BK}
\end{equation}
Moreover,
\begin{equation}
\|e^{it\omega_\alpha(D)}h\|_{\mathcal B_K}
\lesssim
\langle t\rangle^{-3/2}
\sum_{|\rho|\leq2}\|x^\rho h\|_{H^M}.
\label{eq:weighted-linear-BK}
\end{equation}
\end{lemma}
\begin{proof}
Estimate \eqref{eq:WM1-linear-BK} follows immediately by applying
Proposition \ref{prop:dyadic-dispersive} to each dyadic block. Indeed,
the low-frequency contribution is bounded by
\(\langle t\rangle^{-3/2}\|h\|_{L^1}\), while the high-frequency sum
is bounded by
\[
\langle t\rangle^{-3/2}
\sum_{k\geq1}2^{(K+5/2-M)k}\|h\|_{W^{M,1}},
\]
which converges because \(M=K+4\).
It remains to prove the weighted \(L^2\) version
\eqref{eq:weighted-linear-BK}.
For \(P_{\leq0}h\), the compact-frequency estimate in
Proposition \ref{prop:dyadic-dispersive}, together with
\[
\|g\|_{L^1}
\lesssim
\|\langle x\rangle^2g\|_{L^2},
\]
gives
\begin{equation}
\|e^{it\omega_\alpha(D)}P_{\leq0}h\|_{L^\infty}
\lesssim
\langle t\rangle^{-3/2}
\sum_{|\rho|\leq2}\|x^\rho h\|_{L^2}.
\label{eq:low-linear-BK-proof}
\end{equation}
For \(k\geq1\), the high-frequency dyadic estimate gives
\[
2^{Kk}
\|e^{it\omega_\alpha(D)}P_kh\|_{L^\infty}
\lesssim
\langle t\rangle^{-3/2}
2^{(K+5/2)k}\|P_kh\|_{L^1}.
\]
Commuting \(x^\rho\), \(|\rho|\leq2\), through \(P_k\) and applying
Cauchy--Schwarz yields
\[
\|P_kh\|_{L^1}
\lesssim
2^{-Mk}
\sum_{|\rho|\leq2}\|x^\rho h\|_{H^M}.
\]
Since \(M=K+4\), the resulting sum contains
\(\sum_{k\geq1}2^{-3k/2}\), which is convergent. For \(0\leq t\leq1\),
the same conclusion follows from Bernstein and Sobolev embedding.
\end{proof}
Combining Lemma~\ref{lem:weighted-linear-decay} with the weighted
Hardy estimate \eqref{eq:weighted-hardy-position-from-derivative},
we obtain, for the normalized representative of \(h\),
\begin{align}
\|e^{it\omega_\alpha(D)}h\|_{\mathcal B_K}
+\|e^{it\omega_\alpha(D)}\mathcal D_\alpha h\|_{\mathcal B_K}\lesssim
\langle t\rangle^{-3/2}
\sum_{|\rho|\leq3}
\|x^\rho\mathcal D_\alpha h\|_{H^M},
\qquad \alpha\in\mathcal A.
\label{eq:derivative-data-full-linear-decay}
\end{align}
Here the additional spatial weight allows us to control the two
weighted norms of \(h\) from its spatial derivatives.
We now write the normal-form equation directly in terms of the modes
\(U_\alpha\). Recall that
\begin{equation}
(\partial_t-i\omega_\alpha(D))U_\alpha
=\mathcal Q_\alpha(\mathcal U,\mathcal U),
\label{eq:modal-equation-for-normal-form}
\end{equation}
where
\begin{equation}
\widehat{\mathcal Q_\alpha(\mathcal U,\mathcal U)}(\xi)
:=
\sum_{\beta,\gamma\in\mathcal A}
\int_{\mathbb R^3}
m_{\alpha;\beta,\gamma}(\xi,\eta)
\widehat U_\beta(\eta)
\widehat U_\gamma(\xi-\eta)\,d\eta.
\label{eq:physical-quadratic-modal-operator}
\end{equation}
Let \(\mathfrak B_{\alpha;\beta,\gamma}\) be the bilinear
multiplier with symbol \(b_{\alpha;\beta,\gamma}\), and set
\begin{equation}
\mathfrak B_\alpha(\mathcal U,\mathcal U)
:=\sum_{\beta,\gamma\in\mathcal A}
\mathfrak B_{\alpha;\beta,\gamma}(U_\beta,U_\gamma).
\label{eq:physical-normal-form-correction}
\end{equation}
Define the modified modes by subtracting this quadratic correction:
\begin{equation}
G_\alpha
:=U_\alpha-\mathfrak B_\alpha(\mathcal U,\mathcal U).
\label{eq:physical-modified-unknown}
\end{equation}
This agrees with the previously defined modified profile \(F_\alpha\):
\begin{equation}
G_\alpha(t)=e^{it\omega_\alpha(D)}F_\alpha(t).
\label{eq:modified-mode-profile-relation}
\end{equation}
Indeed, since
\(\Phi_{\alpha;\beta,\gamma}
=-\omega_\alpha(\xi)+\omega_\beta(\eta)
+\omega_\gamma(\xi-\eta)\),
\[
e^{it\omega_\alpha(\xi)}
e^{it\Phi_{\alpha;\beta,\gamma}(\xi,\eta)}
\widehat f_\beta(\eta)\widehat f_\gamma(\xi-\eta)
=\widehat U_\beta(\eta)\widehat U_\gamma(\xi-\eta).
\]
Thus the oscillatory quadratic correction in the profile equation
becomes \(\mathfrak B_\alpha(\mathcal U,\mathcal U)\) after
restoring the linear evolution.
The identity \(i\Phi b=m\) gives
\begin{align}
&(\partial_t-i\omega_\alpha(D))
\mathfrak B_\alpha(\mathcal U,\mathcal U)
\notag\\
&\quad=
\mathcal Q_\alpha(\mathcal U,\mathcal U)
+\sum_{\beta,\gamma\in\mathcal A}
\Big[
\mathfrak B_{\alpha;\beta,\gamma}
\bigl(\mathcal Q_\beta(\mathcal U,\mathcal U),U_\gamma\bigr)
+\mathfrak B_{\alpha;\beta,\gamma}
\bigl(U_\beta,\mathcal Q_\gamma(\mathcal U,\mathcal U)\bigr)
\Big].
\label{eq:quadratic-correction-evolution}
\end{align}
Subtracting this identity from
\eqref{eq:modal-equation-for-normal-form} cancels the quadratic
terms and yields
\begin{equation}
(\partial_t-i\omega_\alpha(D))G_\alpha
=\mathcal R_\alpha
(\mathcal U,\mathcal U,\mathcal U),
\label{eq:G-normal-form-generated-equation}
\end{equation}
where
\begin{align}
\mathcal R_\alpha
(\mathcal U,\mathcal U,\mathcal U)
={}&-\sum_{\beta,\gamma\in\mathcal A}
\mathfrak B_{\alpha;\beta,\gamma}
\bigl(\mathcal Q_\beta(\mathcal U,\mathcal U),U_\gamma\bigr)
\notag\\
&-\sum_{\beta,\gamma\in\mathcal A}
\mathfrak B_{\alpha;\beta,\gamma}
\bigl(U_\beta,\mathcal Q_\gamma(\mathcal U,\mathcal U)\bigr).
\label{eq:explicit-normal-form-generated-term}
\end{align}
Each term is cubic because \(\mathcal Q_\beta\) and
\(\mathcal Q_\gamma\) are quadratic in \(\mathcal U\).
Once \(G_\alpha\) has been estimated, we recover \(U_\alpha\) from
\begin{equation}
U_\alpha
=G_\alpha+\mathfrak B_\alpha(\mathcal U,\mathcal U).
\label{eq:U-from-G}
\end{equation}
The norms used below are still the norms of \(\psi,S\) defined
in \eqref{eq:EN-modal}--\eqref{eq:YK-modal}.
The inverse modal formula
\eqref{eq:modal-coordinate-definition} gives
\begin{align}
&\sum_{\alpha\in\mathcal A}
\|\mathcal D_\alpha U_\alpha\|_{H^N}
\lesssim E_N,
\notag\\
&\sum_{\alpha\in\mathcal A}
\|\mathcal D_\alpha U_\alpha\|_{\mathcal B_K}
\lesssim Z_K,
\notag\\
&\sum_{\alpha\in\mathcal A}
\|U_\alpha\|_{\mathcal B_K}
\lesssim Y_K+Z_K.
\label{eq:modal-norms-from-original-variables}
\end{align}
After taking a modal derivative, the inverse formulas involve only
bounded order-zero multipliers applied to
\(\nabla\psi,\nabla S,\psi_t,S_t\).
The undifferentiated modes are controlled by the same formulas
using \(\psi,S,\psi_t,S_t\).
The following lemma gives bounds for the quadratic correction,
the cubic remainder, and the weighted initial data for \(G\).

\begin{lemma}[Quadratic correction and cubic remainder]
\label{lem:normal-form-cubic-estimates}
Assume \(E_N(t)\leq1\). Then
\begin{align}
\sum_\alpha
\Bigl(
\|\mathcal D_\alpha
\mathfrak B_\alpha(\mathcal U,\mathcal U)(t)\|_{\mathcal B_K}
+
\|\mathfrak B_\alpha(\mathcal U,\mathcal U)(t)\|_{\mathcal B_K}
\Bigr)
&\lesssim
\bigl(Y_K(t)+Z_K(t)\bigr)^2,
\label{eq:B-BK-estimate}
\\
\sum_\alpha
\|\mathcal D_\alpha
\mathfrak B_\alpha(\mathcal U,\mathcal U)(t)\|_{H^N}
&\lesssim
E_N(t)\bigl(Y_K(t)+Z_K(t)\bigr),
\label{eq:B-HN-decay-estimate}
\\
\sum_\alpha
\Bigl(
\|\mathcal R_\alpha(t)\|_{W^{M,1}}
+
\|\mathcal D_\alpha\mathcal R_\alpha(t)\|_{W^{M,1}}
\Bigr)
&\lesssim
E_N(t)^2\bigl(Y_K(t)+Z_K(t)\bigr),
\label{eq:R-WM1-estimate}
\\
\sum_\alpha
\Bigl(
\|\mathcal R_\alpha(t)\|_{H^{K+2}}
+
\|\mathcal D_\alpha\mathcal R_\alpha(t)\|_{H^{K+2}}
+
\|\mathcal D_\alpha\mathcal R_\alpha(t)\|_{H^N}
\Bigr)
&\lesssim
E_N(t)Z_K(t)\bigl(Y_K(t)+Z_K(t)\bigr).
\label{eq:R-Sobolev-estimate}
\end{align}
Here
\[
\mathcal R_\alpha(t)
:=
\mathcal R_\alpha
\bigl(\mathcal U(t),\mathcal U(t),\mathcal U(t)\bigr).
\]
Moreover, the initial data satisfy
\begin{equation}
\sum_\alpha
\sum_{|\rho|\leq2}
\Bigl(
\|x^\rho G_\alpha(0)\|_{H^M}
+
\|x^\rho\mathcal D_\alpha G_\alpha(0)\|_{H^M}
\Bigr)
\lesssim
\delta_0+\delta_0^2.
\label{eq:G0-weighted-estimate}
\end{equation}
\end{lemma}

\begin{proof}
All norms below are evaluated at the same time unless the initial
time is explicitly indicated. We use
\eqref{eq:modal-norms-from-original-variables} to express the
modal bounds in terms of \(E_N,Z_K,Y_K\).
The symbol \eqref{eq:full-bilinear-symbol} shows that each
\(\mathcal Q_\beta\) is a finite sum of terms of the form
\[
\langle D\rangle^{-1}(vw),
\]
where \(v,w\) are components of
\(\mathcal D_\delta U_\delta\) and
\(\mathcal D_\kappa U_\kappa\), respectively.
The multipliers
\(\langle D\rangle^{-1}\),
\(\nabla\langle D\rangle^{-1}\), and
\(\lambda_\sigma(D)\langle D\rangle^{-1}\)
are bounded on \(H^s\) and \(\mathcal B_K\).
The tame Sobolev product estimate and the algebra property of
\(\mathcal B_K\) therefore give
\begin{align}
\sum_\beta
\Bigl(
\|\mathcal Q_\beta(\mathcal U,\mathcal U)\|_{H^s}
+\|\mathcal D_\beta
\mathcal Q_\beta(\mathcal U,\mathcal U)\|_{H^s}
\Bigr)
&\lesssim E_NZ_K,
\label{eq:Q-Hs-detailed}
\\
\sum_\beta
\Bigl(
\|\mathcal Q_\beta(\mathcal U,\mathcal U)\|_{\mathcal B_K}
+\|\mathcal D_\beta
\mathcal Q_\beta(\mathcal U,\mathcal U)\|_{\mathcal B_K}
\Bigr)
&\lesssim Z_K^2,
\label{eq:Q-BK-detailed}
\end{align}
for \(K+2\leq s\leq N\).
In the first estimate one differentiated input is measured in
\(H^s\) and the other in \(L^\infty\); in the second, both are
measured in \(\mathcal B_K\).
Applying \eqref{eq:normal-form-pointwise-tame-Hs} and
\eqref{eq:normal-form-pointwise-tame-BK}, together with
\eqref{eq:modal-norms-from-original-variables}, proves
\eqref{eq:B-HN-decay-estimate} and \eqref{eq:B-BK-estimate}.
We next estimate the cubic remainder using
\eqref{eq:explicit-normal-form-generated-term}.
For a typical term
\(\mathfrak B_{\alpha;\beta,\gamma}
(\mathcal Q_\beta,U_\gamma)\), the differentiated tame estimate gives
\begin{align*}
\|\mathcal D_\alpha
\mathfrak B_{\alpha;\beta,\gamma}
(\mathcal Q_\beta,U_\gamma)\|_{H^N}
\lesssim{}&
\|\mathcal D_\beta\mathcal Q_\beta\|_{H^N}
\bigl(\|U_\gamma\|_{\mathcal B_K}
+\|\mathcal D_\gamma U_\gamma\|_{\mathcal B_K}\bigr)
\\
&+
\bigl(\|\mathcal Q_\beta\|_{\mathcal B_K}
+\|\mathcal D_\beta\mathcal Q_\beta\|_{\mathcal B_K}\bigr)
\|\mathcal D_\gamma U_\gamma\|_{H^N}.
\end{align*}
The terms with the two inputs interchanged satisfy the same bound.
Using \eqref{eq:Q-Hs-detailed}--\eqref{eq:Q-BK-detailed} and summing
over the indices, we obtain
\[
\sum_\alpha\|\mathcal D_\alpha\mathcal R_\alpha\|_{H^N}
\lesssim E_NZ_K(Y_K+Z_K)+E_NZ_K^2
\lesssim E_NZ_K(Y_K+Z_K).
\]
This also controls the \(H^{K+2}\) norm of
\(\mathcal D_\alpha\mathcal R_\alpha\).
For the undifferentiated remainder, the explicit decomposition
\eqref{eq:normal-form-operator-decomposition} gives, with \(s=K+2\),
\begin{align*}
\|\mathfrak B_{\alpha;\beta,\gamma}(q,u)\|_{H^s}
\lesssim{}&
\bigl(\|q\|_{H^s}+\|\mathcal D_\beta q\|_{H^s}\bigr)
\bigl(\|u\|_{\mathcal B_K}
+\|\mathcal D_\gamma u\|_{\mathcal B_K}\bigr)
\\
&+
\bigl(\|q\|_{\mathcal B_K}
+\|\mathcal D_\beta q\|_{\mathcal B_K}\bigr)
\|\mathcal D_\gamma u\|_{H^s}.
\end{align*}
To obtain this inequality, split
\(\langle D\rangle u=u+(\langle D\rangle-1)u\), use that
\((\langle D\rangle-1)u\) is controlled by \(\nabla u\), and
apply the tame product estimate. The term with no derivative on
\(u\) is controlled by \(\|q\|_{H^s}\|u\|_{L^\infty}\).
Taking \(q=\mathcal Q_\beta\), \(u=U_\gamma\), and applying
the source estimates proves
\[
\sum_\alpha\|\mathcal R_\alpha\|_{H^{K+2}}
\lesssim E_NZ_K(Y_K+Z_K).
\]
Together these estimates give \eqref{eq:R-Sobolev-estimate}.

For \eqref{eq:R-WM1-estimate}, use the same explicit decomposition
and the representation of \(\mathcal Q_\beta\) as a sum of
\(\langle D\rangle^{-1}(vw)\).
The Bessel potential \(\langle D\rangle^{-1}\) has an integrable
kernel and is bounded on \(L^p\), \(1\leq p\leq\infty\).
At the \(L^1\) endpoint we may use
\[
\|\nabla\langle D\rangle^{-1}H\|_{W^{M,1}}
\lesssim\|H\|_{W^{M+1,1}},
\qquad
\|\langle D\rangle H\|_{W^{M,1}}
\lesssim\|H\|_{W^{M+2,1}}.
\]
The latter follows from
\(\langle D\rangle=(1-\Delta)\langle D\rangle^{-1}\).
Thus at most \(M+2\) additional spatial derivatives need to be
distributed among the three underlying factors, two of which
already carry the modal derivatives supplied by the quadratic source.
The order-zero multipliers in
\((\langle D\rangle-1)u\) are applied to \(\nabla u\) in
\(L^2\) or \(\mathcal B_K\).
In each resulting term, place the two factors with the largest
derivative orders in \(L^2\), and the remaining factor in
\(L^\infty\). If the remaining modal input carries no derivative,
keep it in \(L^\infty\). Every factor placed in \(L^2\) then
carries a spatial or modal derivative and is controlled by \(E_N\).
The factor placed in \(L^\infty\) has at most
\(\lfloor(M+2)/3\rfloor=4\) additional derivatives, and is
controlled by \(Y_K+Z_K\), since \(K=6\).
When a Bessel potential acts on the quadratic product, use its
\(L^1\) boundedness if both inner factors are in \(L^2\), and
its \(L^2\) boundedness if one inner factor is in \(L^\infty\).
H\"older's inequality, with \(M+2\leq N\), now yields
\[
\sum_\alpha
\Bigl(
\|\mathcal R_\alpha\|_{W^{M,1}}
+\|\mathcal D_\alpha\mathcal R_\alpha\|_{W^{M,1}}
\Bigr)
\lesssim E_N^2(Y_K+Z_K).
\]
Finally, we estimate the initial data for \(G\).
The initial-data bounds and the weighted Hardy inequality imply
\begin{equation}
\sum_\alpha\sum_{|\rho|\leq2}
\|x^\rho U_\alpha(0)\|_{H^{M+1}}
\lesssim\delta_0.
\label{eq:weighted-modal-initial-HM1}
\end{equation}
Indeed, Hardy controls the weighted \(H^M\) norms of
\(U_\alpha(0)\), and the weighted \(H^M\) norms of its spatial
derivatives supply the additional derivative.
Commuting \(x^\rho\), \(|\rho|\leq2\), through
\eqref{eq:normal-form-operator-decomposition} differentiates only
the smooth symbols of \(\langle D\rangle\),
\(\langle D\rangle^{-1}\), and their compositions with
\(\mathcal D_\alpha\). Each symbol derivative lowers its order,
and all lower spatial weights are included in the same sum.
The Sobolev product estimate therefore gives
\begin{align*}
&\sum_\alpha\sum_{|\rho|\leq2}
\Bigl(
\|x^\rho\mathfrak B_\alpha
(\mathcal U(0),\mathcal U(0))\|_{H^M}
\\
&\hspace{3cm}
+\|x^\rho\mathcal D_\alpha\mathfrak B_\alpha
(\mathcal U(0),\mathcal U(0))\|_{H^M}
\Bigr)
\\
&\qquad\lesssim
\left(
\sum_\alpha\sum_{|\rho|\leq2}
\|x^\rho U_\alpha(0)\|_{H^{M+1}}
\right)^2
\lesssim\delta_0^2.
\end{align*}
Combining this with
\[
G_\alpha(0)
=U_\alpha(0)
-\mathfrak B_\alpha(\mathcal U(0),\mathcal U(0))
\]
and the corresponding identity after applying \(\mathcal D_\alpha\)
proves \eqref{eq:G0-weighted-estimate}.
\end{proof}
\subsection{Bootstrap and Scattering}
We now prove the main theorem.
\begin{proof}[Proof of Theorem~\ref{thm:small-data-global}]
\medskip
\noindent\textit{Bootstrap improvement.}
Let \([0,T_*)\) be the maximal interval of existence obtained by
iterating Proposition~\ref{prop:local-continuation}.
The initial data estimates give
\[
E_N(0)\leq\delta_0,
\qquad
Y_K(0)+Z_K(0)\leq C\delta_0.
\]
The second bound follows from the weighted Hardy estimate and
Sobolev embedding. Fix a sufficiently large constant \(A\geq1\),
independent of \(\delta_0\). By continuity, the following bounds
hold on some interval \([0,T]\), with \(T<T_*\):
\begin{align}
E_N(t)&\leq2A\delta_0,
\label{eq:bootstrap-EN}
\\
Z_K(t)&\leq2A\delta_0\langle t\rangle^{-3/2},
\label{eq:bootstrap-ZK}
\\
Y_K(t)&\leq2A\delta_0\langle t\rangle^{-3/2}.
\label{eq:bootstrap-YK}
\end{align}
We will replace \(2A\) by \(A\) in these estimates.
We take \(\delta_0\) small enough that \(2A\delta_0\leq1\),
so that Lemma~\ref{lem:normal-form-cubic-estimates} applies
throughout this interval. All constants below are independent
of \(A\), \(\delta_0\), and \(T\).
The energy estimate \eqref{eq:EN-gronwall-ZK} gives
\begin{align}
E_N(t)
&\leq C_E E_N(0)
\exp\left(C\int_0^t Z_K(s)\,ds\right)
\notag\\
&\leq C_E\delta_0\exp(CA\delta_0),
\label{eq:bootstrap-energy-preimprovement}
\end{align}
because \(\int_0^\infty\langle s\rangle^{-3/2}\,ds<\infty\).
Choose \(A\geq2C_E\). If \(\delta_0\) is sufficiently small
that \(\exp(CA\delta_0)\leq2\), then
\begin{equation}
E_N(t)\leq A\delta_0,
\qquad 0\leq t\leq T.
\label{eq:bootstrap-energy-improved}
\end{equation}
We next estimate the modified modes \(G_\alpha\).
Duhamel's formula for
\eqref{eq:G-normal-form-generated-equation} reads
\begin{equation}
G_\alpha(t)
=e^{it\omega_\alpha(D)}G_\alpha(0)
+\int_0^t e^{i(t-s)\omega_\alpha(D)}
\mathcal R_\alpha(s)\,ds.
\label{eq:G-duhamel}
\end{equation}
Apply \eqref{eq:weighted-linear-BK} separately to
\(G_\alpha(0)\) and \(\mathcal D_\alpha G_\alpha(0)\).
The initial estimate \eqref{eq:G0-weighted-estimate} then yields
\begin{align}
\sum_\alpha\Bigl(
&\|e^{it\omega_\alpha(D)}G_\alpha(0)\|_{\mathcal B_K}
+\|e^{it\omega_\alpha(D)}
\mathcal D_\alpha G_\alpha(0)\|_{\mathcal B_K}
\Bigr)
\notag\\
&\lesssim
(\delta_0+\delta_0^2)\langle t\rangle^{-3/2}.
\label{eq:G-linear-part-decay}
\end{align}
Only the two spatial weights supplied by
\eqref{eq:G0-weighted-estimate} are used here.
For the integral term, split the time interval at \(s=t/2\).
On \([0,t/2]\), the elapsed time \(t-s\) is comparable to
\(t\). Using \eqref{eq:WM1-linear-BK} and
\eqref{eq:R-WM1-estimate}, we obtain
\begin{align}
\sum_\alpha\Biggl(
&\left\|\int_0^{t/2}
e^{i(t-s)\omega_\alpha(D)}\mathcal R_\alpha(s)\,ds
\right\|_{\mathcal B_K}
\notag\\
&+\left\|\int_0^{t/2}
e^{i(t-s)\omega_\alpha(D)}
\mathcal D_\alpha\mathcal R_\alpha(s)\,ds
\right\|_{\mathcal B_K}\Biggr)
\notag\\
&\lesssim
\int_0^{t/2}\langle t-s\rangle^{-3/2}
E_N(s)^2\bigl(Y_K(s)+Z_K(s)\bigr)\,ds
\notag\\
&\lesssim
A^3\delta_0^3\langle t\rangle^{-3/2}
\int_0^{t/2}\langle s\rangle^{-3/2}\,ds
\notag\\
&\lesssim A^3\delta_0^3\langle t\rangle^{-3/2}.
\label{eq:G-early-time}
\end{align}
On \([t/2,t]\), use unitarity in \(H^{K+2}\) and the embedding
\(H^{K+2}(\mathbb R^3)\hookrightarrow\mathcal B_K\).
Estimate \eqref{eq:R-Sobolev-estimate} gives
\begin{align}
\sum_\alpha\Biggl(
&\left\|\int_{t/2}^{t}
e^{i(t-s)\omega_\alpha(D)}\mathcal R_\alpha(s)\,ds
\right\|_{\mathcal B_K}
\notag\\
&+\left\|\int_{t/2}^{t}
e^{i(t-s)\omega_\alpha(D)}
\mathcal D_\alpha\mathcal R_\alpha(s)\,ds
\right\|_{\mathcal B_K}\Biggr)
\notag\\
&\lesssim
\int_{t/2}^{t}
E_N(s)Z_K(s)\bigl(Y_K(s)+Z_K(s)\bigr)\,ds
\notag\\
&\lesssim
A^3\delta_0^3\int_{t/2}^{t}\langle s\rangle^{-3}\,ds
\notag\\
&\lesssim A^3\delta_0^3\langle t\rangle^{-2}.
\label{eq:G-late-time}
\end{align}
Both estimates hold for all \(t\geq0\); for \(t\leq2\),
the time integrals are bounded and the factors
\(\langle t\rangle^{-3/2}\) and \(\langle t\rangle^{-2}\)
are comparable to one.
Since \(\mathcal D_\alpha\) commutes with the linear group,
\eqref{eq:G-duhamel} and the preceding estimates imply
\begin{equation}
\sum_\alpha\Bigl(
\|G_\alpha(t)\|_{\mathcal B_K}
+\|\mathcal D_\alpha G_\alpha(t)\|_{\mathcal B_K}
\Bigr)
\lesssim
\bigl(\delta_0+A^3\delta_0^3\bigr)
\langle t\rangle^{-3/2},
\label{eq:G-final-decay}
\end{equation}
where we used \(\delta_0\leq1\).
To recover the norms of \(\psi\) and \(S\), the reconstruction
formulas \eqref{eq:modal-reconstruction-position}--
\eqref{eq:modal-reconstruction-velocity} give
\[
Y_K(t)+Z_K(t)
\lesssim\sum_\alpha\Bigl(
\|U_\alpha(t)\|_{\mathcal B_K}
+\|\mathcal D_\alpha U_\alpha(t)\|_{\mathcal B_K}
\Bigr).
\]
Now use \(U_\alpha=G_\alpha+
\mathfrak B_\alpha(\mathcal U,\mathcal U)\), together with
\eqref{eq:B-BK-estimate}, to obtain
\begin{equation}
Y_K(t)+Z_K(t)
\leq C\bigl(\delta_0+A^3\delta_0^3\bigr)
\langle t\rangle^{-3/2}
+C\bigl(Y_K(t)+Z_K(t)\bigr)^2.
\label{eq:ZK-before-absorption}
\end{equation}
The last term is controlled directly by the bootstrap assumptions:
\[
\bigl(Y_K(t)+Z_K(t)\bigr)^2
\leq16A^2\delta_0^2\langle t\rangle^{-3}
\leq16A^2\delta_0^2\langle t\rangle^{-3/2}.
\]
Consequently,
\begin{equation}
Y_K(t)+Z_K(t)
\leq
\bigl(C_0+C_1A^2\delta_0+C_1A^3\delta_0^2\bigr)
\delta_0\langle t\rangle^{-3/2}.
\label{eq:bootstrap-decay-preimprovement}
\end{equation}
Choose \(A\) large enough to satisfy the earlier requirements
and \(A\geq2C_0\). With this \(A\) fixed, decrease
\(\delta_0\) so that the earlier smallness conditions hold and
\[
C_1A^2\delta_0+C_1A^3\delta_0^2\leq A/2.
\]
Then
\begin{equation}
Y_K(t)+Z_K(t)
\leq A\delta_0\langle t\rangle^{-3/2},
\qquad 0\leq t\leq T.
\label{eq:bootstrap-ZK-improved}
\end{equation}
Since both norms are nonnegative, this improves each of
\eqref{eq:bootstrap-ZK} and \eqref{eq:bootstrap-YK}.
The improved bounds and continuity extend the estimates to
every \(t<T_*\). If \(T_*\) were finite, the uniform bound
\(E_N(t)\leq A\delta_0\) would allow us to restart the local
solution at any \(t_0<T_*\), with a lifespan bounded below
independently of \(t_0\), by
Proposition~\ref{prop:local-continuation}.
Taking \(t_0\) sufficiently close to \(T_*\) would extend the
solution beyond \(T_*\). Hence \(T_*=\infty\).
The bounds \eqref{eq:main-global-bound} and
\eqref{eq:main-derivative-decay} follow, with \(A\) now fixed.

\medskip
\noindent\textit{Scattering.}
The global bounds and \eqref{eq:R-Sobolev-estimate} imply
\begin{equation}
\sum_\alpha\Bigl(
\|\mathcal R_\alpha(t)\|_{H^{K+2}}
+\|\mathcal D_\alpha\mathcal R_\alpha(t)\|_{H^N}
\Bigr)
\lesssim\delta_0^3\langle t\rangle^{-3}.
\label{eq:R-integrable-scattering}
\end{equation}
In particular, both norms are integrable in time.
Recall that \(F_\alpha(t)=e^{-it\omega_\alpha(D)}G_\alpha(t)\).
Equation \eqref{eq:G-normal-form-generated-equation} gives
\begin{equation}
\partial_tF_\alpha(t)
=e^{-it\omega_\alpha(D)}\mathcal R_\alpha(t).
\label{eq:F-scattering-derivative}
\end{equation}
Define the asymptotic profile by
\begin{equation}
f_{\alpha,+}
:=G_\alpha(0)
+\int_0^\infty
e^{-is\omega_\alpha(D)}\mathcal R_\alpha(s)\,ds.
\label{eq:scattering-profile-definition}
\end{equation}
The integral converges in \(H^{K+2}\), and its modal derivative
converges in \(H^N\). The initial term has the same regularity
by \eqref{eq:G0-weighted-estimate},
\eqref{eq:B-HN-decay-estimate}, and the definition of \(G\).
No additional assumption on the initial data is needed.
Subtracting \eqref{eq:scattering-profile-definition} from the
integrated profile equation gives
\[
F_\alpha(t)-f_{\alpha,+}
=-\int_t^\infty
e^{-is\omega_\alpha(D)}\mathcal R_\alpha(s)\,ds.
\]
By unitarity and \eqref{eq:R-integrable-scattering},
\begin{equation}
\sum_\alpha
\|\mathcal D_\alpha(F_\alpha(t)-f_{\alpha,+})\|_{H^N}
\lesssim
\delta_0^3\int_t^\infty\langle s\rangle^{-3}\,ds
\lesssim\delta_0^3\langle t\rangle^{-2}.
\label{eq:F-scattering-rate}
\end{equation}
The quadratic correction also tends to zero in the derivative
energy norm. Indeed, \eqref{eq:B-HN-decay-estimate} gives
\begin{equation}
\sum_\alpha
\|\mathcal D_\alpha(U_\alpha(t)-G_\alpha(t))\|_{H^N}
\lesssim E_N(t)\bigl(Y_K(t)+Z_K(t)\bigr)
\lesssim\delta_0^2\langle t\rangle^{-3/2}.
\label{eq:normal-form-correction-scattering}
\end{equation}
Since \(G_\alpha(t)=e^{it\omega_\alpha(D)}F_\alpha(t)\),
we conclude that
\begin{equation}
\sum_\alpha
\left\|\mathcal D_\alpha\left(
U_\alpha(t)-e^{it\omega_\alpha(D)}f_{\alpha,+}
\right)\right\|_{H^N}
\lesssim
\delta_0^2\langle t\rangle^{-3/2}
+\delta_0^3\langle t\rangle^{-2}.
\label{eq:modal-scattering-rate}
\end{equation}
We now construct the linear solution appearing in the theorem.
Using the reconstruction formulas, set
\begin{align}
\widetilde{\psi}(t)
&:=\sum_{\sigma,\varepsilon\in\{\pm1\}}
e^{it\varepsilon\lambda_\sigma(D)}
f_{(\sigma,\varepsilon),+},
\notag\\
\widetilde{S}(t)
&:=-i\sum_{\sigma,\varepsilon\in\{\pm1\}}
\varepsilon\sigma\,
e^{it\varepsilon\lambda_\sigma(D)}
f_{(\sigma,\varepsilon),+}.
\label{eq:scattering-linear-solution}
\end{align}
Each summand evolves by its corresponding linear mode, so
\((\psi_{\mathrm{lin}},S_{\mathrm{lin}})\) solves
\eqref{eq:linearized-system}.
The conjugacy relations of the original modes pass to
\(f_{\alpha,+}\); hence this linear solution is real-valued.
The position and velocity reconstruction formulas now imply
\begin{align*}
&\|\partial_t(\psi-\widetilde{\psi}(t)\|_{H^N}
+\|\nabla(\psi-\widetilde{\psi}(t)\|_{H^N}
\\
&\quad
+\|\partial_t(S-\widetilde{S}(t)\|_{H^N}
+\|\nabla(S-\widetilde{S}(t)\|_{H^N}
\\
&\lesssim
\delta_0^2\langle t\rangle^{-3/2}
+\delta_0^3\langle t\rangle^{-2}.
\end{align*}
This proves \eqref{eq:physical-scattering} and completes the proof
of Theorem~\ref{thm:small-data-global}.
\end{proof}
\begin{remark}[Use of the spatial weights]
\label{rem:no-XM-needed}
The spatial weights enter the argument through
\eqref{eq:G0-weighted-estimate}, which supplies two weighted
\(H^M\) norms of \(G_\alpha(0)\) and
\(\mathcal D_\alpha G_\alpha(0)\).
These norms control the free evolution.
The cubic integral is estimated using \(W^{M,1}\) on
\([0,t/2]\) and \(H^{K+2}\) on \([t/2,t]\).
Thus the bootstrap closes without estimates for spatially
weighted profiles at positive times.
\end{remark}
\section*{Acknowledgments}
We would like to thank Silu Yin for helpful discussions.

\section*{Funding}
The research of Ben Duan was supported in part by the National Key R\&D Program of China
(No.~2024YFA1013303) and the National Natural Science Foundation of China
(No.~12271205). The research of Rongrong Yan was supported in part by the National Key R\&D Program of China (No.~2024YFA1013303) and by the Scientific Research Foundation for Ph.D. Candidates of the Jilin Provincial Department of Education (No.~JJKH20262211BS). Rongrong Yan also acknowledges support from the China Scholarship Council (CSC, No.~202606170018).

\section*{Conflict of Interest}
All authors declare that they have no conflict of interest.

\section*{Data Availability}
No data was used for the research described in the article.


\begin{thebibliography}{9}
\bibitem{Antonelli2009}
P. Antonelli and P. Marcati, On the finite energy weak solutions to a system in quantum fluid dynamics, Comm. Math. Phys. {\bf 287} (2009), no.~2, 657--686.

\bibitem{Antonelli2012}
P. Antonelli and P. Marcati, The quantum hydrodynamics system in two space dimensions, Arch. Ration. Mech. Anal. {\bf 203} (2012), no.~2, 499--527.

\bibitem{Antonelli2024}
P. Antonelli, P. Marcati and R. Scandone, Existence and stability of almost finite energy weak solutions to the quantum Euler-Maxwell system, J. Math. Pures Appl. (9) {\bf 191} (2024), Paper No. 103629, 38 pp.

\bibitem{Antonelli2017}
P. Antonelli and S. Spirito, Global existence of finite energy weak solutions of quantum Navier-Stokes equations, Arch. Ration. Mech. Anal. {\bf 225} (2017), no.~3, 1161--1199.

\bibitem{AudiardHaspot2017}
C. Audiard and B. Haspot, Global well-posedness of the Euler-Korteweg system for small irrotational data, Comm. Math. Phys. {\bf 351} (2017), no.~1, 201--247.

\bibitem{AudiardHaspot2018}
C. Audiard and B. Haspot, From the Gross-Pitaevskii equation to the Euler Korteweg system, existence of global strong solutions with small irrotational initial data, Ann. Mat. Pura Appl. (4) {\bf 197} (2018), no.~3, 721--760.

\bibitem{BenzoniGavageDanchinDescombes2007}
S. Benzoni-Gavage, R. Danchin and S. Descombes, On the well-posedness for the Euler-Korteweg model in several space dimensions, Indiana Univ. Math. J. {\bf 56} (2007), no.~4, 1499--1579.

\bibitem{Christodoulou1986}
D. Christodoulou, Global solutions of nonlinear hyperbolic equations for small initial data, Comm. Pure Appl. Math. {\bf 39} (1986), no.~2, 267--282.

\bibitem{DuanYan}
B. Duan, J. Li, B. Guo and R. Yan,
A relativistic quantum Euler--Poisson system derived from the
Klein--Gordon--Poisson equation: hyperbolic--elliptic structure,
arXiv:2509.10084, 2025.

\bibitem{GermainMasmoudi2014}
P. Germain and N. Masmoudi,
Global existence for the Euler-Maxwell system,
Ann. Sci. \'Ec. Norm. Sup\'er. (4) {\bf 47} (2014), no.~3, 469--503.

\bibitem{GermainMasmoudiShatah2012}
P. Germain, N. Masmoudi and J. Shatah,
Global solutions for the gravity water waves equation in dimension 3,
Ann. of Math. (2) {\bf 175} (2012), no.~2, 691--754.

\bibitem{GermainMasmoudiShatah2015}
P. Germain, N. Masmoudi and J. Shatah,
Global existence for capillary water waves,
Comm. Pure Appl. Math. {\bf 68} (2015), no.~4, 625--687.

\bibitem{Guo1998}
Y. Guo, Smooth irrotational flows in the large to the Euler-Poisson system in $\bold R^{3+1}$, Comm. Math. Phys. {\bf 195} (1998), no.~2, 249--265.

\bibitem{GuoIonescuPausader2016}
Y. Guo, A.~D. Ionescu and B. Pausader, Global solutions of the Euler-Maxwell two-fluid system in 3D, Ann. of Math. (2) {\bf 183} (2016), no.~2, 377--498.

\bibitem{GuoPausader2011}
Y. Guo and B. Pausader, Global smooth ion dynamics in the Euler-Poisson system, Comm. Math. Phys. {\bf 303} (2011), no.~1, 89--125.

\bibitem{GuoPengWang2008}
Z. Guo, L. Peng and B. Wang, Decay estimates for a class of wave equations, J. Funct. Anal. 254 (2008), no.~6, 1642--1660.

\bibitem{GustafsonNakanishiTsai2009}
S. Gustafson, K. Nakanishi and T.-P. Tsai,
Scattering theory for the Gross--Pitaevskii equation in three dimensions,
Commun. Contemp. Math. {\bf 11} (2009), no.~4, 657--707.

\bibitem{Klainerman1985}
S. Klainerman, Uniform decay estimates and the Lorentz invariance of the classical wave equation, Comm. Pure Appl. Math. {\bf 38} (1985), no.~3, 321--332.

\bibitem{LinWu2012}
C.-K. Lin and K.-C. Wu, Hydrodynamic limits of the nonlinear Klein-Gordon equation, J. Math. Pures Appl. (9) {\bf 98} (2012), no.~3, 328--345.

\bibitem{OzawaTsutayaTsutsumi1995}
T.~Ozawa, K.~Tsutaya, and Y.~Tsutsumi,
Normal form and global solutions for the Klein--Gordon--Zakharov
equations,
Ann. Inst. H. Poincar\'e Anal. Non Lin\'eaire
\textbf{12} (1995), no.~4, 459--503.

\bibitem{Shatah1985}
J.~Shatah,
Normal forms and quadratic nonlinear Klein--Gordon equations,
Comm. Pure Appl. Math. \textbf{38} (1985), 685--696.
\end{thebibliography}
\end{document}